\documentclass[12pt]{amsart}
\usepackage{amscd}

\newtheorem{thm}{Theorem}[section]

\theoremstyle{definition}

\def \ph{\varphi}

\def \refeq#1{equation (\ref{#1})}
\def \ra{\rightarrow}

\def \hom{\mbox{\rm Hom}}
\def \ie{\hbox{\it i.e.}}

\def \tns{\otimes}

\def \k{\mbox{$\mathfrak K$}}

\def \C{\mbox{$\mathbb C$}}
\def \Z{\mbox{$\mathbb Z$}}

\def\br#1#2{\lbrack#1,#2\rbrack}

\def\zt{\mbox{$\Z_2$}}

\def\inv{^{-1}}
\def\d{d}
\def\td{\tilde\d}

\def\A{\mbox{$\mathcal A$}}
\def\B{\mbox{$\mathcal B$}}
\def\L{L}

\def\LA{\mbox{$\L_{\A}$}}
\def\LB{\mbox{$\L_{\B}$}}

\def\m{\mbox{$\mathfrak m$}}
\def\a{\mbox{$\mathfrak a$}}

\def\coder{\operatorname{Coder}}

\def\linf{\mbox{$L_\infty$}}
\def\and{\mbox{ \rm and }}

\def\s#1{(-1)^{#1}}
\DeclareMathOperator*{\invlim}{\overleftarrow{\rm lim}}

\def\htns{\hat\tns}

\def\htns{\hat\tns}

\def\phd#1#2{\ph^{#1}_{#2}}
\def\two{w_2}
\def\one{w_1}
\def\tre{w_3}
\def\three{\tre}
\def\thmref#1{Theorem (\ref{#1})}
\def\inv{^{-1}}
\def\psd#1#2{\psi^{#1}_{#2}}
\def\dstar{\mbox{$d_*$}}

\def\dsharp{\mbox{$d_{\sharp}$}}

\def\dt#1{\mbox{$d_{#1}$}}

\author{Derek Bodin}
\address{University of Wisconsin\\
Eau Claire, WI 54702-4004}
\email{bodinda@uwec.edu}
\author{Alice Fialowski}
\address{E\"otv\"os Lor\'and University\\
Budapest, Hungary} \email{fialowsk@cs.elte.hu}
\author{Michael Penkava}
\address{University of Wisconsin\\
Eau Claire, WI 54702-4004}
\email{penkavmr@uwec.edu}
\thanks{The authors were unable to attend the meeting in Seville in
person, but wished to contribute this paper. The research was
 partially supported by
the grants MTA-OTKA-NSF 38453, OTKA T043641 and T043034 and by grants from the
 University of Wisconsin-Eau Claire}
\title[\linf\ algebras of dimension $2|1$]{Classification and versal
deformations of \linf\ algebras on a $2|1$-dimensional space}
\begin{document}
\setlength{\multlinegap}{0pt}
\begin{abstract}
This article explores $\Z_2$-graded
\linf\ algebra structures on a $2|1$-dimensional
vector space. The reader should note that our convention
on the parities is the opposite of the usual one, because we define our
structures on the symmetric coalgebra of the parity reversion of a
space, so our $2|1$-dimensional \linf\ algebras correspond to the usual
$1|2$-dimensional algebras.

We give a complete classification of all structures with
a nonzero degree 1 term. We also classify all degree 2 codifferentials, which
is the same as a classification of all $1|2$-dimensional \zt-graded Lie
algebras. For each of these algebra structures, we calculate the cohomology
and a miniversal deformation.

\end{abstract}
\date\today
\maketitle

\section{Introduction}

\linf\ algebras - or strongly homotopy Lie algebras - were first described
in \cite{ss} and have recently been the focus of much attention both in
mathematics (\cite{sta3}, \cite{hsch}) and in mathematical physics
(\cite{ls}, \cite{z}, \cite{aksz}, \cite{bfls}, \cite{RW}, \cite{M},
\cite{fls}).

In the physics literature, one usually considers $\Z_2$-graded spaces.
That's the case in our consideration also: throughout this paper, all spaces
will be $\Z_2$-graded, and we will work in the
parity reversed definition of the \linf\ structure.
In \cite{fp2}, all \linf\ algebras of
dimension less than or equal to 2 were classified, in \cite{fp3} miniversal
 deformations for all \linf\ structures on a
space of dimension $0|3$ ($0$ even and $3$ odd dimension)- which
correspond to ordinary Lie algebras - were constructed, and in
\cite{fp4}, all \linf\ algebras of dimension $1|2$ were classified.

The picture in the $2|1$ dimensional case is
more complicated than  $1|2$-dimensional algebras,
because the space of $n$-cochains on a $1|2$
dimensional space has dimension $6|6$ for $n>1$, while the space
of $n$-cochains on a $2|1$- dimensional space has dimension
$3n+2|3n+1$, making it more difficult to classify the
nonequivalent structures.  Accordingly, we give a complete
classification here of only those \linf\ algebras which are extensions
of degree 1 coderivations, which are, as it turns out, equivalent to
degree 1 coderivations, and those which are extensions of degree 2 coderivations,
in other words,  extensions of \zt-graded Lie algebras as \linf\
algebras. For each of these algebras we either construct a miniversal
deformation, or give the necessary steps for the construction.

\section{Basic Definitions}

\subsection{\linf\ algebras}

We work in the framework of the parity reversion $W=\Pi V$ of the usual
vector space $V$ on which an \linf\ algebra structure is defined,
because in the $W$ framework,  an \linf\ structure is simply an odd
coderivation $d$ of the symmetric coalgebra $S(W)$, satisfying
$d^2=0$,  in other words, it is an odd codifferential in the \zt-graded
Lie algebra of coderivations of $S(W)$.  As a consequence, when
studying \zt-graded Lie algebra structures on $V$, the parity is
reversed,  so that a $2|1$-dimensional vector space $W$ corresponds to
a $1|2$-dimensional \zt-graded Lie structure on $V$.  Moreover,  the
\zt-graded anti-symmetry of the Lie bracket on $V$ becomes the
\zt-graded symmetry of the associated coderivation $d$ on $S(W)$.

A formal power series $d=d_1+\cdots$, with $d_i\in L_i=\hom(S^i(W),W)$ determines
an element in $L=\hom(S(W),W)$, which is naturally identified with
$\coder(S(W))$, the space of coderivations of the symmetric
coalgebra $S(W)$. Thus $L$ is a \zt-graded Lie algebra.
An odd element $d$ in $L$ is called a  \emph{codifferential} if
$\br dd=0$. We also say that $d$ is
an \linf\ structure on $W$.

If $g=g_1+\cdots\in\hom(S(W),W)$, and $g_1:W\ra W$ is invertible,
then $g$ determines a coalgebra automorphism of $S(W)$ in a natural
way, which we will denote by the same letter $g$. Moreover, every
coalgebra automorphism is determined in this manner. Two codifferentials
$d$ and $d'$ are said to be \emph{equivalent} if there is a coalgebra
automorphism $g$ such that $d'=g^*(d)=g\inv d g$.

A  detailed description of \linf\ algebras can be obtained in
\cite{lm,ls,pen3,pen4}. The study of examples of \linf\ algebra structures
in \cite{fp2,fp3,fp4} may be useful to the reader,  but we
intend this article to be as self contained as possible.

\subsection{Equivalent codifferentials and extensions}

We will construct all equivalence classes of \linf\ structures on $W$.
We shall use the following facts, which
are established in \cite{fp1,bfp2}, to aid in the classification.

If $d$ is an \linf\ structure on $W$, and $d_N$ is the first
nonvanishing term in $d$, then $d_N$ is itself a codifferential, which we
call the \emph{leading term of $d$}, and we say that $d$ is
\emph{an extension} of $d_N$. Define the cohomology operator $D$
by $D(\ph)=\br\ph {d_N}$, for $\ph\in L$. Then the following
formula holds for any extension $d$ of $d_N$ as an \linf\
structure, and all $n\ge N$:
\begin{equation}\label{extension}
D(d_{n+1})=-\frac12\sum_{k=N+1}^{n}\br{d_k}{d_{n+N-k+1}}.
\end{equation}
Note that the terms on the right all have index less than $n+1$.
If a coderivation $d$ has been constructed up to terms of degree
$m$, satisfying \refeq{extension} for $n=1\dots m-1$, then the
right hand side of \refeq{extension} for $n=m$ is automatically a
cocycle. Thus $d$ can be extended to the next level precisely when
the cocycle given by the right hand side is trivial. There may be many
nonequivalent extensions, because the term $d_{n+1}$ which we add
to extend the coderivation is only determined up to a cocycle.  An
extension $d$ of $d_N$ is given by any coderivation whose leading term is
$d_N$, which satisfies
\refeq{extension} for every $m=N+1\dots$. The theory here is
parallel to the theory of formal deformations of an algebra structure;
the extension of a codifferential $d_N$ to a more complicated
codifferential $d$ resembles the process of extending an infinitesimal
deformation to a formal one.

Classifying the extensions of $d_N$ can be quite complicated.
However, the following theorem often makes it easy to classify the
extensions.
\begin{thm}\label{zerocohomology}
If the cohomology $H^n(d_N)=0$, for $n>N$, then any extension of
$d_N$ to a \linf\ structure $d$ is equivalent to the structure
$d_N$.
\end{thm}

Before classifying the extensions of a codifferential $d_N$, we
need to classify the codifferentials in $L_N$ up to equivalence,
that is, we need to study the \emph{moduli space} of degree $N$
codifferentials. A \emph{linear automorphism} of $S(W)$ is an
automorphism determined by an isomorphism $g_1:W\ra W$. If $g$ is
an arbitrary automorphism, determined by maps $g_i:S^i(W)\ra W$,
and $W$ is finite dimensional, then $g_1$ is an isomorphism, so
this term alone induces an automorphism of $S(W)$ which we call
the \emph{linear part} of $g$.

The following theorem simplifies the classification of equivalence
classes of codifferentials in $L_N$.
\begin{thm}\label{simpleauto}
If $d$ and $d'$ are two codifferentials in $L_N$, and $g$ is an
equivalence between them, then the linear part of $g$ is also an
equivalence between them.
\end{thm}
Thus we can restrict ourself to linear automorphisms when
determining the equivalence classes of elements in $L_N$.

We will also use the following result.
\begin{thm}\label{simpleequiv}
Suppose that $d$ and $d'$ are equivalent codifferentials.  Then
their leading terms have the same degree and are equivalent.
\end{thm}

As a consequence of these theorems, we proceed to classify the
codifferentials as follows. First,  find all equivalence classes
of codifferentials of degree $N$. For each equivalence class,
study the equivalence classes of extensions of the codifferential.

Let us first establish some basic notation for the cochains.
Suppose $W=\langle w_1,w_2,w_3\rangle$, with $w_1$ an odd element and $w_2,w_3$
 even elements.
If $I=\{i_1,i_2,i_3\}$ is a
multi-index, with $i_1$ and $i_2$ either zero or one, let
$w_I=w_1^{i_1}w_2^{i_2}w_3^{i_3}$. For simplicity, we will denote
$w_I$ by $I$.
Then for $n\ge 1$,
\begin{align*}
(S^{n}(W))_e=&\langle (0,p,n-p)|0\le p\le n\rangle,\qquad |(S^{n}(W))_e|=n+1\\
(S^{n}(W))_o=&\langle (1,q,n-q-1)|0\le q\le n-1\rangle,\qquad |(S^{n}(W))_0|=n
\end{align*}

If $\lambda$ is a
linear automorphism of $S(W)$, then in terms of the standard basis of $W$, its restriction to $W$ has matrix
\begin{equation}\label{linauto}
\lambda=\begin{pmatrix}
q&0&0\\
0&r&t\\
0&s&u\\
\end{pmatrix}
\end{equation}
where $q(ru-st)\ne 0$.
We will sometimes express $\lambda$ by the submatrix $\bigl(\begin{smallmatrix}r&t\\s&u\end{smallmatrix}\bigr)$.
It is useful to note that for a linear automorphism
$$\lambda(w^I)=\lambda(w_1)^{i_1}\lambda(w_2)^{i_2}\lambda(w_3)^{i_3},$$
so that
\begin{equation}\label{lambdaxy}
\lambda(1,x,y)=\sum_{i=0}^x\sum_{j=0}^y(1,i+j,x+y-i-j)\binom xi \binom yj qr^is^{x-i}t^ju^{y-j}.
\end{equation}

Let $L_n:=\hom(S^n(W),W)$. Define
\begin{equation*}
\ph^I_j(w_J)=I!\delta^I_J w_j,
\end{equation*}
where $I!=i_1!i_2!i_3!$. If we let  $|I|=i_1+i_2+i_3$, then
$L_n=\langle \ph^I_j, |I|=n \rangle$. If $\ph$ is odd, we denote
it by the symbol $\psi$ to make it easier to distinguish the even
and odd elements. Then
\begin{align*}
(L_{n})_e&=\langle\ph^{1,q,n-q-1}_1,\ph^{0,p,n-p}_2,\ph^{0,p,n-p}_3|1\le q\le n-1,1\le p\le n\rangle\\
(L_{n})_o&=\langle\psi^{1,q,n-q-1}_2,\psi^{1,q,n-q-1}_3,\psi^{0,p,n-p}_1|1\le q\le n-1,1\le p\le n\rangle,
\end{align*}
so that $|L_n|=3n+2|3n+1$.

\subsection{ Versal Deformations}

For the classical formal deformation theory we refer to \cite{Ger1}.
Versal deformation theory was first worked out for the case of Lie
algebras in \cite{fi,ff2} and then extended to \linf\ algebras in
\cite{fp1}.

An augmented local ring $\A$ with maximal ideal $\m$ will be called an
\emph{infinitesimal base} if  $\m^2=0$, and a  \emph{formal base} if $\A=\invlim_n \A/\m^n$.
A deformation of an \linf\ algebra structure $d$ on $W$
with  base given by a local ring $\A$ with augmentation $\epsilon:\A\ra\k$, where
$\k$ is the field over which $W$ is defined, is an $\A$-\linf\ structure $\td$ on $W\htns\A$
such that the morphism of $\A$-\linf\ algebras
$\epsilon_*=1\tns\epsilon:\LA=L\tns \A\ra L\tns\k=L$ satisfies
$\epsilon_*(\td)=d$.  (Here $W\htns\A$ is an appropriate completion of $W\tns\A$.)
The deformation is called infinitesimal (formal) if $\A$ is an infinitesimal (formal) base.

In general, the cohomology $H(D)$ of $d$ given by the operator
$D:L\ra L$ with $D(\ph)=\br\ph d$ may not be finite dimensional.
However, $L$ has a natural
 filtration $L^n=\prod_{i=n}^\infty L_i$,
which induces a filtration $H^n$ on the cohomology, because
 $D$
respects the filtration. Then $H(D)$ is of finite type if $H^n/H^{n+1}$
is finite dimensional.
 Since this is always true when $W$ is finite
dimensional, the examples we study here will always
 be of finite
type.  A set $\delta_i$ will be called a \emph{basis of the cohomology},
if any element $\delta$
 of the cohomology can be expressed uniquely
as a formal sum $\delta=\delta_i a^i$. If we identify
 $H(D)$ with a
subspace of the space of cocycles $Z(D)$, and we choose a \emph{basis}
$\beta_i$ of the
 coboundary space $B(D)$, then any element
$\zeta\in\Z(D)$ can be expressed uniquely as a sum
$\zeta=\delta_i a^i +\beta_i b^i$.

For each $\delta_i$, let $u^i$ be a parameter of opposite parity.  Then the infinitesimal
deformation $d^1=d+\delta_i u^i$, with base $\A=\k[u_i]/(u_iu_j)$ is universal in the sense that
if $\td$ is any infinitesimal deformation with base $\B$, then there is a unique homomorphism
$f:\A\ra\B$, such that the morphism $f_*=1\tns f:\LA\ra\LB$ satisfies $f_*(\td)\sim d$.

For formal deformations, there is no universal object in the sense above. A
\emph{versal deformation}
is a deformation $d^\infty$ with formal base $\A$ such that if $\td$ is any formal deformation
with base $\B$, then there is some morphism $f:\A\ra\B$ such that $f_*(d^\infty)\sim\td$. If $f$
is unique whenever $\B$ is infinitesimal, then the versal deformation is called \emph{miniversal}.
In \cite{fp1}, we constructed a miniversal deformation for \linf\ algebras with finite type
cohomology.

The method of construction is as follows.  Define a coboundary
operator $D$ by $D(\ph)=[\ph,d]$. First, one constructs the universal
infinitesimal deformation $d^1=d+\delta_i u^i$, where $\delta_i$
is a graded basis of the cohomology $H(D)$ of $d$, or more
correctly, a basis of a subspace of the cocycles which projects
isomorphically to a basis in cohomology, and $u^i$ is a parameter
whose parity is opposite to $\delta_i$.  The infinitesimal
assumption that the products of parameters are equal to zero gives
the property that $[d^1,d^1]=0$. Actually, we can express
\begin{equation*}
[d^1,d^1]=\s{\delta_j(\delta_i+1)}[\delta_i,\delta_j]u^iu^j=\delta_k
a^k_{ij}u^iu^j +\beta_k b^k_{ij}u^iu^j,
\end{equation*}
where $\beta_i$ is a basis of the coboundaries, because the
bracket of $d^1$ with itself is a cocycle. Note that the right
hand side is of degree 2 in the parameters, so it is zero up to
order 1 in the parameters.

If we suppose that $D(\gamma_i)=-\frac12\beta_i$, then by
replacing $d^1$ with $$d^2=d^1+\gamma_kb^k_{ij}u^iu^j,
$$ one obtains
\begin{equation*}
[d^2,d^2]=\delta_k a^k_{ij}u^iu^j+ 2[\delta_l
u^l,\gamma_kb^k_{ij}u^iu^j]+[\gamma_kb^k_{ij}u^iu^j,\gamma_lb^l_{ij}u^iu^j]
\end{equation*}
Thus we are able to get rid of terms of degree 2 in the
coboundary terms $\beta_i$, but those which involve the cohomology
terms $\delta_i$ can not be eliminated.  This gives rise to a set
of second order relations on the parameters. One continues this
process, taking the bracket of the \emph{$n$-th order deformation}
$d^n$, adding some higher order terms to cancel coboundaries,
obtaining higher order relations, which extend the second order
relations.

Either the process continues indefinitely, in which
case the miniversal deformation is expressed as a formal power
series in the parameters, or after a finite number of steps, the
right hand side of the bracket is zero after applying the $n$-th
order relations. In this case, the miniversal deformation is
simply the $n$-th order deformation. In any case, we obtain a set
of relations $R_i$ on the parameters, one for each $\delta_i$, and
the algebra $A=\C[[u^i]]/(R_i)$ is called the base of the
miniversal deformation. Examples of the construction of miniversal
deformations can be found in \cite{ff3,ff2,fp2,fp3}.
\newpage

\section{Classification of Codifferentials}
Let us compute the brackets of all odd cochains with each other.
\begin{align*}
\br{\psi^{0,p,n-p}_1}{\psi^{1,q,m-q-1}_2}&=\ph^{1,p+q-1,n-p+m-q-1}_1p+\ph^{0,p+q,n-p+m-q-1}_2\\
\br{\psi^{0,p,n-p}_1}{\psi^{1,q,m-q-1}_3}&=\ph^{1,p+q,n-p+m-q-2}_1(n-p)+\ph^{0,p+q,n-p+m-q-1}_3\\
\br{\psi^{0,p,n-p}_1}{\psi^{0,q,m-q}_1}&=0\qquad\qquad
\br{\psi^{1,p,n-p-1}_2}{\psi^{1,q,m-q-1}_2}=0\\
\br{\psi^{1,p,n-p-1}_2}{\psi^{1,q,m-q-1}_3}&=0\qquad\qquad
\br{\psi^{1,p,n-p-1}_3}{\psi^{1,q,m-q-1}_3}=0
\end{align*}
Suppose that
\begin{equation}
d=\sum_{p=0}^n\psi^{0,p,n-p}_1a_p+\sum_{q=0}^{n-1}\psi^{1,q,n-q}_2b_q+\psi^{1,q,n-q}_3c_q,
\end{equation}
where we sum over all odd codifferentials of degree $n$.
Then using the above, we compute that
\begin{multline}
\br dd=\
\sum_{p=0}^{n}\sum_{q=0}^{n-1}
\ph^{1,p+q-1,2n-p-q-1}_1pa_pb_q+\ph^{0,p+q,2n-p-q-1}_2a_pb_q
\\+
\ph^{1,p+q,2n-p-q-2}_1(n-p)a_pc_q+\ph^{0,p+q,2n-p-q-1}_3a_pc_q
\end{multline}
We claim that either all coefficients $a_p$ must vanish or all coefficients $b_q$ and $c_q$ must vanish.
For if $p$ and $q$ are the least indices  for which $a_p$ and $b_q$ do not vanish,  then
there is only one term  in the sum above of the form $\ph^{0,p+q,2n-p-q-1}_2$, which would be a contradiction
because its coefficient $a_pb_q$ must vanish.

As a consequence of this observation,  we note that codifferentials of degree $n$ fall into two distinct
families,  those of the
\emph{first kind}
\begin{equation}
\sum_{q=0}^{n-1}\psi^{1,q,n-q}_2b_q+\psi^{1,q,n-q}_3c_q.
\end{equation}
and those of the \emph{second kind}
\begin{equation}
d=\sum_{p=0}^n\psi^{0,p,n-p}_1a_p.
\end{equation}

Moreover, any expression of either kind gives a codifferential.  Thus we have determined all codifferentials
of degree $N$.  However, the process of classification requires that we determine the equivalence classes of
codifferentials under the action of the automorphism group of the symmetric coalgebra, and we are a long
way away from this classification at this stage.

Several things can be said in general.  First, let us suppose that $d$ is of the second kind. Then from the brackets
computed already,  we note that the odd $d$-cocycles are precisely the odd cochains of the second kind.
The space of odd cocycles
has dimension $n+1$,  which means that the space of even $d$-coboundaries has dimension $2n$.  Also, if $\ph$ is any
even cocycle, then its bracket with $d$ is an odd cocycle of the second kind.  Precise computation of the cohomology
depends on solving a linear system of equations whose coefficients depend on the coefficients in $d$.

Similarly, if $d$ is of the first kind, then the odd $d$-cocycles are the ones of the first kind, and thus the dimension
of the space of odd cocycles is $2n$, and the dimension of the space of even coboundaries is $n+1$.  The bracket
of any even cocycle with $d$ is a cocycle of the first kind.  The cohomology can be computed by solving a system of
linear equations in coefficients depending on the coefficients of $d$.

Since there are $2n$ coefficients in a codifferential of the first kind, and $n+1$ coefficients in a codifferential
of the second kind,  there are potentially a lot of equivalence classes of codifferentials.  The main strategy involved
in classification is to reduce the number of independent variables to a manageable number.

It is useful to compute the brackets of even and odd cochains.
\begin{align*}
[\phd{1,p,m-p}1,\psd{0,q,n-q}1]=&\psd{0,p+q,m+n-p-q}1\\
[\phd{0,p,m-p}2,\psd{0,q,n-q}1]=&-\psd{0,p+q-1,m+n-p-q}1q\\
[\phd{0,p,m-p}3,\psd{0,q,n-q}1]=&-\psd{0,p+q,m+n-p-q-1}1(n-q)\\
[\phd{1,p,m-p}1,\psd{1,q,n-q}2]=&-\psd{1,p+q,m+n-p-q}2\\
[\phd{0,p,m-p}2,\psd{1,q,n-q}2]=&\psd{1,p+q-1,m+n-p-q}2(p-q)\\
[\phd{0,p,m-p}3,\psd{1,q,n-q}2]=&-\psd{1,p+q,m+n-p-q-1}2(n-q)+\psd{1,p+q-1,m+n-p-q}3p\\
[\phd{1,p,m-p}1,\psd{1,q,n-q}3]=&-\psd{1,p+q,m+n-p-q}3\\
[\phd{0,p,m-p}2,\psd{1,q,n-q}3]=&\psd{1,p+q,m+n-p-q-1}2(m-p)-\psd{1,p+q-1,m+n-p-q}3q\\
[\phd{0,p,m-p}3,\psd{1,q,n-q}3]=&\psd{1,p+q,m+n-p-q-1}3(m-p-(n-q))
\end{align*}
Notice that bracket of any even cochain with an odd cochain of a certain type is an odd cochain of the same
type. This is very important in what follows, because this fact means that there is no mixing of types occurring
in the cohomology of a codifferential of a fixed type.

Later on,  in the computation of miniversal deformations,  we will also need to know the brackets of even cochains with
each other, so we include these calculations now.
\begin{align*}
[\phd{1,p,m-p}1,\phd{1,q,n-q}1]=&0\\
[\phd{0,p,m-p}2,\phd{1,q,n-q}1]=&-\phd{1,p+q-1,m+n-p-q}1q\\
[\phd{0,p,m-p}3,\phd{1,q,n-q}1]=&-\phd{1,p+q,m+n-p-q-1}1(n-q)\\
[\phd{0,p,m-p}2,\phd{0,q,n-q}2]=&\phd{0,p+q-1,m+n-p-q}2(p-q)\\
[\phd{0,p,m-p}3,\phd{0,q,n-q}2]=&-\phd{0,p+q,m+n-p-q-1}2(n-q)+\phd{0,p+q-1,m+n-p-q}3p\\
[\phd{0,p,m-p}3,\phd{0,q,n-q}3]=&\phd{0,p+q,m+n-p-q-1}3(m-p-(n-q))
\end{align*}

Let us call the degree of the leading term of a codifferential the
\emph{order} of that codifferential.  We begin with
a classification of all codifferentials of order 1.

\section{Classification and miniversal deformations of codifferentials with $d_1\ne0$}

Let us suppose that $d$ is an odd, degree 1 codifferential of the first kind. Then
$d=\psi^{1,0,0}_2a_1 +\psi^{1,0,0}_3a_2$ for some constants $a_1$ and $a_2$. To see that $d$ is equivalent to
$d'=\psi^{1,0,0}_2$, let $t$ and $u$ be such that $a_1t+a_2u\ne0$.  Suppose that
$g=\begin{pmatrix}a_1&t\\a_2&u\end{pmatrix}$.  Then
\begin{align*}
dg(\one)=d(\one)=\two a_1 +\tre a_2=g(\two)=gd'(\one).
\end{align*}
Since $dg(\two)=gd'(\two)=0$ and $dg(\tre)=gd'(\tre)=0$, it follows that $d'$ and $d$ are
equivalent. Thus every codifferential of the first kind is equivalent to $\psi^{1,0,0}_2$.

Now let us study the cohomology of the codifferential $d=\psi^{1,0,0}_2$.  We define the coboundary operator
$D$ by $D(\ph)=\br{\ph}d$. Then
Computing brackets, we see that
\begin{align*}
D(\ph^{0,p,n-p}_2)&=\psi^{1,p-1,n-p}_2p,\quad
D(\psi^{0,p,n-p}_1)=\ph^{1,p-1,n-p}_1p+\ph^{0,p,n-p}_2\\
D(\ph^{0,p,n-p}_3)&=\psi^{1,p-1,n-p}_3p,\quad
D(\psi^{1,q,n-q-1}_2)=0\\
D(\ph^{1,q,n-q-1}_1)&=-\psi^{1,q,n-q-1}_2,\quad
D(\psi^{1,q,n-q-1}_3)=0
\end{align*}
Note that for the $n$-cochains above,  $p$ ranges from $0$ to $n$, while $q$ ranges only from $0$ to $n-1$.
It is easy to see that $\psi^{1,q,n-q-1}_2$ and $\psi^{1,q,n-q-1}_3$ give a basis of the odd cocycles,
and since both of these types are evidently coboundaries,  of $\ph^{0,p,n-p}_2$ and $\ph^{0,p,n-p}_3$, resp.,
where $q=p-1$, all odd cocycles are coboundaries. Similarly, if we let $q=p-1$, then every even cocycle
is a linear combination of elements of the form $\ph^{1,p-1,n-p}_1p+\ph^{0,p,n-p}_2$, and since these elements
are coboundaries,  it follows that all even cocycles are coboundaries.  Thus the cohomology of $d$ is zero,
and we know by \thmref{zerocohomology} that all extensions of $d$ are equivalent to $d$. This completes
the picture for codifferentials of the first kind of degree 1.

If $d$ is a codifferential of the second kind of degree 1, it is of the form $d=\psi^{0,1,0}_1a_1+\psi^{0,0,1}_1a_2$.
We show that it is equivalent to $d'=\psi^{0,1,0}_1$.  For suppose that $b_1$ and $b_2$ are chosen so that
$a_1b_1+a_2b_2=1$. Then if $g$ is given by $\begin{pmatrix}b_1&-a_2\\b_2&a_1\end{pmatrix}$, we have
\begin{align*}
dg(\two)&=d(\two b_1+\tre b_2)=\one (a_1b_1+a_2b_2)=\one=gd'(\two)\\
dg(\tre)&=d(-\two a_2 +\tre a_1)=-a_2a_1+a_1a_2=0=gd'(\tre)
\end{align*}
Now,  we study the cohomology induced by $d=\psi^{1,0,0}_2$.
Calculating coboundaries, we have
\begin{align*}
D(\ph^{0,p,n-p}_2)=&-\psi^{0,p,n-p}_1,&
D(\psi^{0,p,n-p}_1)=&0\\
D(\ph^{0,p,n-p}_3)=&0,&
D(\psi^{1,q,n-q-1}_2)=&\ph^{1,q,n-q-1}_1+\ph^{0,q+1,n-q-1}_2\\
D(\ph^{1,q,n-q-1}_1)=&\psi^{1,q+1,n-q-1}_1,&
D(\psi^{1,q,n-q-1}_3)=&\ph^{0,q+1,n-q-1}_3
\end{align*}
It is not difficult to see from this table that the cohomology of this
codifferential is also zero. Thus every extension of a codifferential
of the second kind is equivalent to the original codifferential.

The picture for codifferentials of degree 1 is very simple. First, the
classification into equivalence classes is easy, and then,  since the
cohomology vanishes, the classification of extensions is immediate.
There are exactly two equivalence classes of codifferentials of order
1. Because the cohomology vanishes,  there are no nontrivial deformations
of the infinity algebra structures,  so the miniversal deformation of the
\linf\ algebras determined by
degree 1 codifferentials coincides with the codifferentials.


We will study codifferentials of degree two next, with an eventual aim
to classify all extensions of such codifferentials to \linf\ algebra
structures. The classification of extensions will be deferred to a future
work, and we will confine ourselves here to a study of the versal deformations
of the \linf\ algebra structures determined by degree 2 (quadratic) codifferentials.
Since quadratic codifferentials of the symmetric algebra precisely correspond
to the \zt-graded Lie algebra structures on the parity reversion of our space $W$,
what we are really doing here is giving a complete classification of \zt-graded
Lie algebras on a $1|2$-dimensional vector space, and studying their versal deformations
as \linf\ algebras.  The versal deformations of these \linf\ algebras as superalgebras
are given by considering only the deformations induced by quadratic cocycles,
so are an easily identifiable part of the versal deformation we will study.

\section{Codifferentials of Degree 2 of the First Kind}

Let us suppose that $d=\psd{1,1,0}2x+\psd{1,1,0}3a+\psd{1,0,1}2b
+\psd{1,0,1}3c$, and let us call the multi-index $(x,a,b,c)$ the type
of the codifferential. Let us say that a codifferential is of type
$(x,a,b,c)$ whenever it is equivalent to  a codifferential of that
type, so that the type of a codifferential is not unique. Our goal is
to show that the equivalence classes of codifferentials reduce to only
a few simple types.  Let us first remark that if we express $d$ as a
matrix of the form $d=\begin{pmatrix}x&b\\a&c\end{pmatrix}$, then if
$d'=g\inv d g$, then its matrix is simply the product of the matrices
expressing $g\inv$, $d$ and $g$, multiplied by the scalar $q$, where
$g(\one)=\one q$.

If $x\ne0$ then by applying a simple coalgebra automorphism, one can
assume that it is equal to one. Similarly, if $a\ne0$ one can assume it
is also 1. Thus if both $x$ and $a$ are non zero, our codifferential is
of type $(1,1,b,c)$. We will show later that we can express
codifferentials of this type in an even simpler form, but first we
examine what possibilities have not been covered by our considerations.

If $x=0$, but $c\ne 0$, then by interchanging the roles of $\two$ and
$\three$ one can replace it with an equivalent one whose $\psd{1,1,0}2$
has nonzero coefficient. Similarly, if $x=0$ or $a=0$ and both $b$ and
$c$ do not vanish, then by the same interchange, we can see that the
codifferential has type $(1,1,b,c)$ as well. This observation leads to
the following possibilities,

If $x\ne 0$ but $a=0$, then we can assume that either $b=0$ or $c=0$
(or both). This gives the possible types $(1,0,0,c)$ or $(1,0,b,0)$. If
$b\ne 0$, then by a simple transformation, the type $(1,0,b,0)$ can be
reduced to type $(1,0,1,0)$.

The only other types which could arise have both the $x$ and $c$
coefficients vanishing, so they are of type $(0,a,b,0)$. If both $a$
and $b$ do not vanish, they can be adjusted so we obtain type
$(0,1,1,0)$, and if one of the two vanishes but the other does not, we
obtain type $(0,1,0,0)$.

Actually, this myriad of types can be much reduced as we shall see
shortly. Let us examine the type $(1,1,b,c)$ and show that in most
cases it can be reduced to type $(1,0,0,c)$.

Let $d$ be of type $(1,1,b,c)$. Then if $d'=g\inv d g$, we compute
\begin{equation}\label{genform}
d'=
\begin{pmatrix}
\dfrac{q(ru+bus-rt-cst)}{ru-ts}&\dfrac{q(ut+bu^2-t^2-ctu)}{ru-ts}\\
\\
-\dfrac{q(sr+bs^2-r^2-crs)}{ru-ts}&
-\dfrac{q(ts+bus-rt-cru)}{ru-ts}
\end{pmatrix}
\end{equation}
Now, either $r$ and $u$ both do not vanish, or $s$ and $t$ both do not.  Let us assume the former,
and put $x=s/r$ and $y=t/u$.  Substituting in the matrix for $d'$, we obtain
\begin{equation*}
d'=
\begin{pmatrix}
\dfrac{q(1+bx-y-cxy)}{1-xy}&\dfrac{q(y+b-y^2-cy)\frac ur}{1-xy}\\\\
-\dfrac{q(x+bx^2-1-cx)\frac ru}{1-xy}&
-\dfrac{q(xy+bx-y-c)}{1-xy}
\end{pmatrix}
\end{equation*}
Our goal is to remove the off diagonal terms without violating the
condition $xy\ne 1$.  The terms vanish precisely when the equations
$y+b-y^2-cy$ and $x+bx^2-1-cx$ are both equal to zero. When $b\ne0$,
these equations are quadratic in $y$ and $x$ respectively, with
solutions
$$x_{\pm}=\frac{c-1\pm\sqrt{(1-c)^2+4b}}{2b}\qquad y_{\pm}=\frac{1-c\pm\sqrt{(1-c)^2+4b}}2.
$$
Oddly enough, we compute $x_+y_+=x_-y_-=1$, which is just what we want
to avoid.  On the other hand, $x_+y_-=1$, if and only if
$b=-\frac{(1-c)^2}4$. Assuming otherwise, we can eliminate the off
diagonal terms, so that after applying a simple automorphism,  we can
reduce it to type $(1,0,0,c')$, where  $c'$ is given by some rational
expression in $b$ and $c$.

On the other hand, when $b=0$, then the quadratic in $x$ reduces to a
linear expression, which is zero when $x=\frac1{1-c}$. Of course, $x$
is not well defined if $c=1$, so let us first assume otherwise. Now
$y=0$ is a solution of the quadratic equality for $y$, and substituting
the expressions for $x$ and $y$ into the first and fourth terms yields
that our codifferential is equivalent to one of type $(1,0,0,c)$ where
the $c$ in this expression is the same as the $c$ occurring in the type
$(1,1,0,c)$. In fact, it is also clear that even when $b\ne 0$, one can
reduce any expression of type $(1,1,b,c)$ to the type $(1,1,0,c')$ by
choosing $y=y_{+}$, and $x=0$, except in the special case when
$b=-\frac{(1-c)^2}4$.

Note that the case $b=0$ and $c=1$ is a special case of the equality
$b=-\frac{(1-c)^2}4$, so all we have left is to consider the case where
this equality holds. Then we certainly can set $y=y_+=\frac{1-c}2$, and
the condition $1-xy\ne 0$ reduces to the inequality $\frac{c-1}2x+1\ne
0$.  If we choose an arbitrary $x$ so this inequality is satisfied,
then it is easy to see that the first and the fourth coefficients of
$d'$ become, simply $q\frac{c+1}2$. This means that when $c\ne-1$, we
can choose $q$ to make the first and fourth coefficients of $d'$ equal
to 1, the third coefficient equal to 0, and by choosing $r/u$
appropriately, the second coefficient equal to 1 as well. Thus we
obtain an element of type $(1,1,0,1)$ unless $c=-1$. One can check that
in this case, we obtain $b=-1$, so the element has type $(1,1,-1,-1)$,
and also it is obvious from this argument that in this case $d$ is
equivalent to the codifferential $d=\psd{1,1,0}3$, so that in
particular $\psd{1,1,0}3$ has type $(1,1,-1,-1)$. It is also easy to
show that both types $(1,1,-1,-1)$ and $(1,1,0,1)$ can never be reduced
to type $(1,0,0,c)$.

We now proceed to show that the other special types, $(1,0,1,0)$, $(0,1,1,0)$ also can be reduced to
type $(1,0,0,c)$.

Type $(1,0,1,0)$ is the same as type $(1,0,0,0)$ and type $(1,1,0,0)$.
To see this, apply the generic linear transformation $g$ to produce
$d'$ as before, and we obtain
\begin{equation*}
d'=\begin{pmatrix}
\dfrac{q(r+s)u}{ru-ts}&\dfrac{q(t+u)u}{ru-ts}\\\\
-\dfrac{q(r+s)s}{ru-ts}&-\dfrac{q(t+u)s}{ru-ts}
\end{pmatrix}
\end{equation*}
If we choose $q=1/2$, $s=t=r=-1$ and $u=1$, we obtain type $(1,1,0,0)$,
and if instead we choose $q=r=u=1$, $s=0$, and $t=-1$, we obtain type
$(1,0,0,0)$.

Type $(0,1,1,0)$ is the same as  type $(1,0,0,-1)$ and type $(1,1,0,-1)$.
To see this, apply the generic linear transformation $g$ and we obtain
\begin{equation*}
d'=
\begin{pmatrix}
\dfrac{q(su-rt)}{ru-ts}&\dfrac{q(u^2-t^2)}{ru-ts}\\\\
-\dfrac{q(s^2-r^2)}{ru-ts}&-\dfrac{q(su-rt)}{ru-ts}
\end{pmatrix}
\end{equation*}
Choose $t=-1$ and $q=r=u=s=1$. Then this becomes $d'=\psd{1,1,0}2-\psd{1,0,1}3$ as desired.  On the other hand,
if $q=s=u=1$, $r=0$,  and $t=-1$, then we obtain
$d'=\psd{1,1,0}2+\psd{1,1,0}3-\psd{1,0,1}3$.

This completes the classification of types of codifferentials.  We have
one family $(1,0,0,c)$ and two special cases, $(1,1,0,1)$ and
$(0,1,0,0)$ which cannot be reduced to elements of this family.

We show that an element of type $(1,0,0,c)$ is equivalent to one of
type $(1,0,0,c')$ precisely when $c'=c^{\pm 1}$, so that the set of
equivalence classes of codifferentials has a one parameter subfamily,
parameterized by the unit disc in $\mathbb C$.  To see this,  apply the
generic linear transformation and we obtain
\begin{equation*}
d'= \begin{pmatrix}
\dfrac{q(ru-cst)}{ru-ts}&-\dfrac{qtu(-1+c)}{ru-ts}\\\\
-\dfrac{qrs(-1+c)}{ru-ts}& \dfrac{q(-ts+rcu)}{ru-ts} \end{pmatrix}
\end{equation*}
It is interesting to note that when $c=1$, the middle two terms drop
out and thus $\psd{1,1,0}2+\psd{1,0,1}3$ is not equivalent to any
codifferential of type $(1,1,b,c)$. Otherwise, if $c\ne0$, let $r=u=0$
and $s=t=1$ and $q=-1/c$ and we obtain type $(1,0,0,1/c)$. To see that
this is the only other type that could occur, note that to cancel the
middle terms, we must have either $r=u=0$ or $s=t=0$, so the claim is
obvious.

For later purposes let us label the codifferentials representing the
equivalence classes of degree 2 codifferentials of the first kind as
follows.
\begin{align*}
&\dstar=\psd{1,1,0}3\\
&\dsharp=\psd{1,1,0}2+\psd{1,1,0}3+\psd{1,0,1}3\\
&\dt{c}=\psd{1,1,0}2+\psd{1,0,1}3c
\end{align*}

\section{Cohomology of Codifferentials of degree 2 of the First Kind}
For a degree 2 codifferential $d$, with cohomology operator
$D=\br{\bullet}d$, the dimension of the cohomology is given by
$h_n=z_n-b_{n-1}$.

\subsection{Cohomology of $\dstar=\psd{1,1,0}3$}\label{0100}
The coboundaries of basic
cochains for \dstar\ are as follows:
\begin{align*}
&D(\phd{1, q, n-q-1}1) =-\psd{1, 1+q, n-q-1}3\\
&D(\phd{0, p, n-p}2) =\psd{1, 1+p, n-p-1}2(n-p)-\psd{1,p, n-p}3\\
&D(\phd{0, p, n-p}3) =\psd{1, 1+p, n-p-1}3(n-p)\\
&D(\psd{0, p, n-p}1) = \phd{1, 1+p, n-p-1}1(n-p)+\phd{0, 1+p, n-p}3\\
&D(\psd{1, q, n-q-1}2) =
D(\psd{1, q, n-q-1}3) =0
\end{align*}
From this table, we see that  $\psd{1, q, n-q-1}2$, $\psd{1, q, n-q-1}3$
are cocycles for $q=0\dots n-1$, and
$\phd{0, p, n-p}3+\phd{1, p, n-p-1}1(n-p)$ is a cocycle for $p=0\dots n$.
Also, $\phd{0, n, 0}2+\phd{0, n-1, 1}3$ is a cocycle, so
$z_n=n+2|2n$,  which means that $b_{n}=n+1|2n$. It follows that $h_n=2|2$ for $n>1$ and $h_1=3|2$.
Moreover,
\begin{align*}
&H^1=\langle \psd{1,0,0}2,\psd{1,0,0}3,\phd{0,1,0}3,\phd{0,0,1}3+\phd{1,0,0}1,
\phd{0,1,0}2+\phd{0,0,1}3\rangle\\
&H^n=\langle \psd{1,0,n-1}2,\psd{1,0,n-1}3,\phd{0,n,0}2+\phd{0,n-1,1}3,
\phd{0,0,n}3+\phd{1,0,n-1}1n\rangle,\quad n>1.
\end{align*}

Because the odd part of the cohomology of \dstar\ does not vanish for
$n>2$,  there are nontrivial extensions of \dstar. We will discuss these extensions
in a later paper.

\subsection{Cohomology of $\dsharp=\psd{1,1,0}2+\psd{1,1,0}3+\psd{1,0,1}3$}

The coboundaries for \dsharp\ are given by
\begin{align*}
D(\phd{1, q, n-q-1}1) &=-\psd{1, 1+q, n-q-1}2-\psd{1, 1+q, n-q-1}3-\psd{1, q, n-q}3\\
D(\phd{0, p, n-p}2) &=\psd{1, 1+p, n-p-1}2(n-p)-\psd{1,p, n-p}2(n-1)-\psd{1,p, n-p}3\\
D(\phd{0, p, n-p}3) &=\psd{1, 1+p, n-p-1}3(n-p)-\psd{1,p, n-p}3(n-1)\\
D(\psd{0, p, n-p}1) &=\phd{1,p,n-p}1n+\phd{1, 1+p, n-p-1}1(n-p)+\phd{0, 1+p, n-p}2\\&+\phd{0,p,n-p+1}2+\phd{0, 1+p, n-p}3\\
D(\psd{1, q, n-q-1}2) &=
D(\psd{1, q, n-q-1}3) =0
\end{align*}
We already know that $\psd{1,q,n-q-1}2$ and $\psd{1,q,n-q-1}3$ give a
basis of the $2n$ odd cocycles. First note that $\ph_3^{0,1,0}$ and
$\ph_2^{0,1,0}+\ph_3^{0,0,1}$ are a basis of the even 1-cocycles. Thus
$z_1=h_1=2|2$ and  $b_1=2|3$.

For $n>1$ it is easy to see that images of the cochains of the form
$\phd{0,p,n-p}3$ are a basis of  the $n+1$-dimensional subspace of
cochains of the form $\psd{1,p,n-p}3$. If we consider the subspace $X$
spanned by elements of the form $\psd{1, 1+q, n-q-1}2$ and
$\psd{1,p,n-p}3$, then $D$ maps the $(2n+1)$-dimensional space spanned
by elements of the form $\phd{1,q,n-q-1}1$ and $\phd{0,p,n-p}3$
bijectively onto $X$.

If $p>0$, it is clear that the image of $\phd{0,p,n-p}2$ lies in $X$.
Thus we obtain a cocycle as a sum of the element $\phd{0,p,n-p}2$ and a
unique linear combination of the elements $\phd{0,p,n-p}3$ and
$\phd{1,q,n-q-1}1$.  Thus there are $n$ independent cocycles generated
by these elements.

For $p=0$, the image of $\phd{0,0,n}2$ does not lie in $X$, so it
cannot contribute to any cocycle.  Thus we see that there are exactly
$n$ independent even cocycles. Thus $z_n=n|2n$ and $b_n=n+1|2n+2$.

This means that $h_n=z_n-b_{n-1}=0$ if $n>2$.  Furthermore, $h_2=2|4
-2|3=0|1$. It is easy to see that $\psd{1,0,1}2$ can be taken as the
basis for $H^2$. Thus we have
\begin{align*}
&H^1=\langle\psd{1,0,0}2,\psd{1,0,0}3, \phd{0,1,0}3,\phd{0,1,0}2+\phd{0,0,1}3\rangle\\
&H^2=\langle \psd{1,0,1}2\rangle\\
&H^n=0,\qquad\text{if $n>2$}
\end{align*}
Because $H^n=0$ for $n>2$, there are no nontrivial extensions of $\dsharp$. We will discuss
the versal deformation of this codifferential in the next section.
\subsection{Cohomology of $\dt c=\psd{1,1,0}2+\psd{1,0,1}3c$}
Since $\dt c$ is equivalent to $\dt{1/c}$, we can assume that $c$ lies
in the unit circle. Thus we will assume that $|c|\le 1$ in the
following. The coboundaries are given by
\begin{align*}
D(\phd{1, q, n-q-1}1) &=-\psd{1, 1+q, n-q-1}2-\psd{1, q, n-q}3c\\
D(\phd{0, p, n-p}2) &=\psd{1,p, n-p}2(p-1+c(n-p))\\
D(\phd{0, p, n-p}3) &=\psd{1,p, n-p}3(p+c(n-p-1))\\
D(\psd{0, p, n-p}1) &=\phd{1,p,n-p}1(p+c(n-p))+\phd{0, 1+p, n-p}2+\phd{0, p, n-p+1}3c\\
D(\psd{1, q, n-q-1}2) &=0\\
D(\psd{1, q, n-q-1}3) &=0\\
\end{align*}
Let $Q_p=p+c(n-p-1)$.  When $Q_p\ne0$, then
\begin{equation*}
\xi_p=\ph_2^{0,p+1,n-p-1}+\ph_3^{0,p,n-p}c+\ph_1^{1,p,n-p-1}Q_p,\quad p=0\dots n-1,
\end{equation*}
give $n$ independent even cocycles, which are obviously coboundaries.
In most cases, the $\xi_p$ give a basis of the even cocycles. However,
when $Q_p=0$ we get an additional even cocycle $\phd{0,p,n-p}3$ which
is never a coboundary, and $\psd{1,p,n-p}3$ is no longer a coboundary.
When this happens the even part of $h_n$ increases by one, and the odd
part of $h_{n+1}$ also increases by one. Note that if $n>1$, for most
values of $c$ it never happens that $Q_p=0$. In fact, if $c$ is not a
nonpositive rational number, then $Q_p$ is never zero when $n>1$.

There is another source of possible even cocycles, given by the terms
$\phd{0,0,n}2$, which is a cocycle if $nc=1$, and $\phd{0,n,0}3$, which
is a cocycle if $n=c$. When this happens, the even part of $z_n$
increases by 1, so the odd part of $b_n$ decreases by 1. Thus we see
again that the even part of $h_n$ and the odd part of $h_{n+1}$ both
increase by 1. Moreover, if $nc=1$, then we need to add $\psd{1,0,n}2$
to the basis of $H^{n+1}$ and if $n=c$, we need to add $\psd{1,n,0}3$
to the basis. If $c$ or its reciprocal is not a positive integer, then
neither of these two cases hold.

\subsubsection{Cohomology for generic values of $c$}

Let us say that $c$ is generic if it is not a nonpositive rational
number, nor is it or its reciprocal a positive integer. If $c$ is
generic,  and $n>1$, then $\xi_p$ are the only even cocycles, so that
$z_n=n|2n$ and thus we have $b_n=n+1|2n+2$. It follows that $h_n=0|0$
for $n>2$. For $n=1$, we always have $Q_0=0$, so we have two even
cocycles, $\phd{0,1,0}2$ and $\phd{0,0,1}3$, which along with the two
odd cocycles $\psd{1,0,0}2$ and $\psd{1,0,0}3$ generically form a basis
for $H^1$. Thus $z_1=h_1=2|2$ in the generic case.

Also, in the generic case,  we obtain $b_1=2|3$.  Since $z_2=2|4$, it
follows that $h_2=0|1$. In fact, $\ph_3^{1,0,1}$ can be taken as a
basis for $H^2$. What this says is that you can deform $\dt c$ in the
direction of the family.

Thus we conclude that for generic values of $c$ we have
\begin{align*}
&H^1=\langle \psd{1,0,0}2,\psd{1,0,0}3,\phd{0,1,0}2,\phd{0,0,1}3\rangle\\
&H^2=\langle \psd{1,0,1}3\rangle\\
&H^n=0,\qquad\text{if $n>2$}
\end{align*}

Note that for generic values of $c$, the picture is the same as for $\dsharp$,
there are no nontrivial extensions, and the pattern for versal deformations will
be seen to be similar to $\dsharp$ as well.

\subsubsection{Cohomology for the special value $c=1$}

In this case both $nc=1$ and $n=c$ hold for $n=1$. Thus we obtain two
additional 1-cohomology classes, given by $\phd{0,0,1}2$ and
$\phd{0,1,0}3$ and $h_1=z_1=4|2$. Thus $b_1=2|1$, so
$h_2=z_2-b_1=2|4-2|1=0|3$. Thus for $c=1$ we have
\begin{align*}
&H^1=\langle \psd{1,0,0}2,\psd{1,0,0}3,\phd{0,1,0}2,\phd{0,0,1}3,\phd{0,0,1}2,\phd{0,1,0}3\rangle\\
&H^2=\langle \psd{1,0,1}3,\psd{1,0,1}2,\psd{1,1,0}3\rangle\\
&H^n=0,\qquad\text{if $n>2$}
\end{align*}
This suggests that somehow there are additional directions in which the
codifferential can be deformed, and we will comment on this later.


\subsubsection{Cohomology when $1/c\ne 1$ is a positive integer}

Let $m=1/c$. Then we have
\begin{align*}
&H^1=\langle \psd{1,0,0}2,\psd{1,0,0}3,\phd{0,1,0}2,\phd{0,0,1}3\rangle\\
&H^2=\langle \psd{1,0,1}3\rangle\\
&H^m=\langle \phd{0,0,m}2\rangle\\
&H^{m+1}=\langle \psd{1,0,m}2\rangle\\
&H^n=0,\qquad\text{otherwise}
\end{align*}
except when $c=1/2$, in which case, since $m=2$, $H^2=\langle
\psd{1,0,1}3,\phd{0,0,2}2\rangle$.

In this case, we see that there are nontrivial extensions of $\dt c$, and the versal deformation
picture is more complicated as well.

\subsubsection{Cohomology when $c=0$}

This case is special because $Q_0=0$ for all $n$. Thus we always have
the even cohomology class $\phd{0,0,n}3$, and the odd cohomology class
$\psd{1,0,n-1}3$. Since $Q_0$ is zero when $n=1$ in all cases, $H^1$ is
not changed from the generic pattern.  Thus
\begin{align*}
&H^1=\langle \psd{1,0,0}2,\psd{1,0,0}3,\phd{0,1,0}2,\phd{0,0,1}3\rangle\\
&H^n=\langle \psd{1,0,n-1}3,\phd{0,0,n}3\rangle,\qquad\text{if $n>1$}
\end{align*}

For this special case, both the extension and the versal deformation pictures are more involved.

\subsubsection{Cohomology when $c$ is a negative rational number}

Let us rewrite the equality $Q_p=0$ in the form $p=\frac{(n-1)c}{c-1}$
When $-1\le c<0$ is rational, note that $0<c/(c-1)\le \tfrac12$, so
that $Q=0$ has an integral solution for $p$ with $0<p<n-1$ for
infinitely many values of $n$. In fact, suppose that $\frac
c{c-1}=\frac rs$, expressed as a fraction in lowest terms.  Then
$\frac{(n-1)c}{c-1}$ is a positive integer precisely when $n=ks+1$ for
a positive integer $k$, in which case we have $kr=\frac{(n-1)c}{c-1}$,
and $0<kr<n-1$.  From this, we can calculate the table of cohomology of
$\dt c$ as follows.
\begin{align*}
&H^1=\langle \psd{1,0,0}2,\psd{1,0,0}3,\phd{0,1,0}2,\phd{0,0,1}3\rangle\\
&H^2=\langle \psd{1,0,1}3\rangle\\
&H^{ks+1}=\langle \phd{0,kr,k(s-r)+1}3\rangle\\
&H^{ks+2}=\langle \psd{1,kr,k(s-r)+1}3\rangle\\
&H^n=0,\qquad\text{otherwise}
\end{align*}
There are nontrivial extensions of this type of codifferential.

\subsubsection{The Moduli Space of Codifferentials of the First Kind}

Let us consider now the deformations of these various types of
codifferentials of degree 2 only as graded Lie algebras, \ie, consider
only $H^2$.  Consider the following table of codifferentials and bases
of the odd part of the second cohomology group:

\begin{equation}\label{table1}
\notag
\begin{array}{ccccc}
\\
\text{Type}&(0,1,0,0)&(1,1,0,1)&(1,0,0,1)&(1,0,0,c)\\
\\
d&\dstar&\dsharp&$\dt1$&$\dt c$\\
\\
(H^2)_o&\psd{1,0,1}2,\psd{1,0,1}3&\psd{1,0,1}2&\psd{1,0,1}2,\psd{1,0,1}3,\psd{1,1,0}3&\psd{1,0,1}3
\end{array}
\end{equation}

There are three special cases, and the generic pattern.  Note that even
though the dimension of $H^2$ is not generic for $\dt 0$, the extra
dimension is even, so does not contribute to the deformations over
$\mathbb C$.

Clearly, there is only one family of codifferentials, so what is going
on with the extra degrees of freedom in the cohomology?  To understand
this better, let us examine the moduli space of codifferentials of
degree 2 in some more detail. Here we use the term moduli space in the
following sense.  The space of all codifferentials of degree 2 is a
variety in a 4 dimensional complex space, preserved under the action of
the group of linear automorphisms of the symmetric coalgebra. A
quotient space of a variety by such a group action is called a moduli
space. The structure of such moduli spaces can be very strange, from a
topological point of view.

Let us parameterize our moduli space by types, and note that since type
$(1,0,0,c)$ is the same as type $(1,0,0,1/c)$, it is natural to think
of the moduli space as the unit disc in $\mathbb C$, with an
identification of the upper semicircle with the bottom.  Then every
point except 1 and -1 have neighborhoods which are discs, but 1 and -1
are orbifold points of degree 2. Of course, we are really describing
the action of the group generated by the transformations $\{z\ra z,
z\ra 1/z\}$ on the Riemann sphere, and identifying our standard points
of the moduli space with the resulting images.

We should like to have some notion of neighborhood of a point in our
moduli space, and the natural notion is to consider two elements of the
moduli space to be close if they have inverse images which are close in
the space of codifferentials. Of course, since any codifferential is
equivalent to any multiple of itself, this would make all
codifferentials close, so we have to be a bit more careful in our
definition.

Consider the standard representatives of the equivalence classes of
codifferentials, which are either $(1,0,0,c)$, $(1,1,0,1)$ and
$(0,1,0,0)$. Let $P$ and $Q$ be equivalence classes.  Let us say that
$Q$ is $\epsilon$ close to $P$ if $Q$ is among the types which occur by
adding coordinates to the standard representation of $P$ of absolute
value no larger than $\epsilon$. Then $P$ is said to be infinitesimally
close to $Q$ if $Q$ is epsilon close to $P$ for all positive values of
$\epsilon$.

For most of our points, the notion of neighborhood we have just
described yields no surprises.  For any standard point $P$ of type
$(1,0,0,c)$ with $c\ne 1$, $\epsilon$ neighborhoods of $P$ for small
values of $\epsilon$ correspond to standard points $(1,0,0,c')$ with
$c'$ close to $c$.

However, for $c=1$, things are quite different. As we have described
before, the one parameter family $(1,1,b,c)$ with $b=-\frac14(c-1)^2$,
contains types $(1,0,0,1)$ and $(0,1,0,0)$ for two special values of
$c$, but gives type $(1,1,0,1)$ otherwise.  It follows that $\dt1$ and
$\dstar$ are infinitesimally close to $\dsharp$. One can check that for
$\epsilon$ small enough, a neighborhood of $\dt1$ contains only this
extra point, along with the points one would usually expect. It is hard
to reconcile the fact that for the codifferential $\dt1$, the dimension
of the cohomology is 3. One might expect one extra dimension for the
deformation in the $\dsharp$ direction, but two extra dimensions are
obtained instead.

Recall that type $(1,1,0,a)$ is the same as type $(1,0,0,a)$ when
$a\ne1$.  This means that a neighborhood of $\dsharp$ looks just like a
neighborhood of $\dt1)$ (minus point $\dt1$).  Note that although
$\dt1$ is infinitesimally close to $\dsharp$, the converse is not true.
Notice that the dimension of the cohomology for $\dsharp$ is just 1,
corresponding to the fact that any small deformation of this
codifferential just gives an ordinary element in the main family.

Finally, consider type $(0,1,0,0)$. Note that $(0,1,\epsilon,0)$ is the
same as type $(1,0,0,-1)$, so $\dstar$ is infinitesimally close to
$\dt{-1}$. It is easy to see that type $(0,1,\epsilon_1,\epsilon_2,0)$
is the same as type $(1,1,\epsilon_1/\epsilon_2^2,0)$, if
$\epsilon_2\ne0$.  One also sees that type $(1,0,0,c)$ is the same as
type $(1,1,-\frac{c}{(c+1)^2},0)$, if $c\ne \pm 1$. When $c=1$, we
obtain type $(1,1,-\frac14,0)$ which is the same as type $(1,1,0,1)$.
Thus $\dstar$ is infinitesimally close to every element of the moduli
space except $\dt1$. Note that the cohomology has odd dimension 2, and
the type $(0,1,\epsilon_1,\epsilon_2)$ corresponds to adding a small
cocycle to $\dstar$.

\section{Miniversal deformations of degree 2 codifferentials of the
first kind}

The most important part of the construction of a miniversal deformation
is the computation of the relations on its base, because they determine
the answer to the classical question: ``Given an infinitesimal
deformation, when does it extend to a formal deformation?".  Sometimes
it is possible to calculate the relations on the base of a miniversal
deformation, without explicitly computing the miniversal deformation.
In most of the examples here, we give explicit computations of the
versal deformation, but for the structure $\dstar$, even the
computation of the second order deformation is quite involved, so a
general formula would be difficult to develop.

\subsection{A miniversal deformation of $\dsharp$}
We study this one first because it is very simple in comparison.
The universal infinitesimal deformation is given by
\begin{align*}
d^1=\dsharp+\psd{1,0,0}2t_1+\psd{1,0,0}3t_2+\phd{0,1,0}3\theta_1+(\phd{0,1,0}2+\phd{0,0,1}3)\theta_2
+\psd{1,0,1}2 t_3.
\end{align*}
To compute the versal deformation,  we first compute
\begin{align*}
\tfrac12[d^1,d^1]=&
-\psd{1,0,0}3(t_1\theta_1+t_2\theta_2)-\psd{1,0,0}2t_1\theta_2+(\psd{1,1,0}2-\psd{1,0,1}3)t_3\theta_1.
\end{align*}
The first two terms are cohomology classes, so give rise to the second order relations
\begin{equation*}
t_1\theta_1+t_2\theta_2=0,\qquad t_1\theta_2=0.
\end{equation*}
The third term is a coboundary, in fact $D(\phd{0,0,1}2)=\psd{1,1,0}2-\psd{1,0,1}3$.
Thus the second order deformation of $\dsharp$ is given by
\begin{equation*}
d^2=d^1-\phd{0,0,1}2t_2\theta_1.
\end{equation*}
Continuing the process, we obtain
\begin{equation*}
\tfrac12[d^2,d^2]=\psd{1,0,0}2t_2t_3\theta_1.
\end{equation*}
As a consequence of the fact that no coboundary terms appear in this bracket,  we see that $d^2$ is a
miniversal deformation of $\dsharp$,  and the relations become
\begin{equation*}
t_1\theta_1+t_2\theta_2=0,\qquad -t_1\theta_2+t_2t_3\theta_1=0.
\end{equation*}
The base $\A$ of the versal deformation is thus
\begin{equation*}
\A=\k[[t_1,t_2,t_3,\theta_1,\theta_2]]/(t_1\theta_1+t_2\theta_2,-t_1\theta_2+t_2t_3\theta_1).
\end{equation*}

\subsection{A miniversal deformation of $\dt c$ for generic values of $c$}

The universal infinitesimal deformation is given by
\begin{equation*}
d^1=\psd{1,1,0}2+\psd{1,0,1}3c+
\psd{1,0,0}2 t^1+\psd{1,0,0}3 t^2 +\phd{0,1,0}2\theta^1+\phd{0,0,1}3\theta^2+\psd{1,0,1}3 t^3.
\end{equation*}
Then
\begin{equation*}
\tfrac12[d^1,d^1]=-\psd{1,0,0}2 t^1\theta^1-\psd{1,0,0}3 t^2\theta^2.
\end{equation*}
Since both of the terms in the bracket are cohomology classes, $d^1$ is already a miniversal deformation
of $\dt c$, and the relations on the base are simply
\begin{equation*}
t^1\theta^1=0,\qquad t^2\theta^2=0.
\end{equation*}

\subsection{A miniversal deformation of $\dt 1$}

The universal infinitesimal deformation is given by
\begin{align*}
d^1=&\psd{1,1,0}2+\psd{1,0,1}3+
\psd{1,0,0}2 t^1+\psd{1,0,0}3 t^2 +\phd{0,1,0}2\theta^1+\phd{0,0,1}3\theta^2+
\phd{0,0,1}2\theta^3\\&+\phd{0,1,0}3\theta^4+\psd{1,0,1}3 t^3+\psd{1,0,1}2 t^4+\psd{1,1,0}3 t^5.
\end{align*}
Then
\begin{align*}
\tfrac12[d^1,d^1]=&-\psd{1,0,0}2 (t^1\theta^1 +t^2\theta^3)-\psd{1,0,0}3 (t^2\theta^2+t^1\theta^4)
\\&-
\psd{1,0,1}2(t^3\theta^3+t^4\theta^1-t^4\theta^2)+\psd{1,1,0}3(t^3\theta^4+t^5\theta^1-t^5\theta^2)
\\&+(\psd{1,1,0}2-\psd{1,0,1}3)(t^4\theta^4- t^5\theta^3).
\end{align*}
This is a rather interesting situation,  because
\begin{equation*}(\psd{1,1,0}2-\psd{1,0,1}3)=D(\phd{1,0,0}1)+2\psd{1,0,1}3,
\end{equation*}
so it is a sum of a coboundary term and a cohomology class.  Therefore,  a second order relation
is $t^4\theta^4- t^5\theta^3=0$, and after taking the second order relations into account, one obtains that
the bracket vanishes identically.  Thus, the versal deformation of $\dt 1$ is still given by the
universal infinitesimal deformation, as in the generic case, but with the relations
\begin{align*}
t^1\theta^1 +t^2\theta^3=
t^2\theta^2+t^1\theta^4=
t^4\theta^4- t^5\theta^3=0\\
t^3\theta^3+t^4\theta^1-t^4\theta^2=
t^3\theta^4+t^5\theta^1-t^5\theta^2=0\\
\end{align*}

\subsection{Deformations of $\dt c$ when $m=1/c\ne 1$ is a positive integer}

The universal infinitesimal deformation is given by \begin{align*}
d^1=&\psd{1,1,0}2+\psd{1,0,1}3c+ \psd{1,0,0}2 t^1+\psd{1,0,0}3 t^2
+\phd{0,1,0}2\theta^1+\phd{0,0,1}3\theta^2+\psd{1,0,1}3 t^3
\\&+\phd{0,0,m}2\theta^3+\psd{1,0,m}2 t^4. \end{align*} Then
\begin{align*} \tfrac12[d^1,d^1]=&-\psd{1,0,0}2
t^1\theta^1-\psd{1,0,0}3 t^2\theta^2
+\phd{0,0,m}2(\theta^1\theta^3-m\theta^2\theta^3)\\
&-\psd{1,0,m}2(mt^3\theta^3+t^4\theta^1-mt^4\theta^2)
-\psd{1,0,m-1}2mt^2\theta^3. \end{align*} The only coboundary term
appearing is $-\psd{1,0,m-1}2=D(\phd{0,0,m-1}2m)$, so the second order
relations are \begin{equation*} t^1\theta^1= t^2\theta^2=
\theta^1\theta^3-m\theta^2\theta^3=
mt^3\theta^3+t^4\theta^1-mt^4\theta^2=0, \end{equation*} and the second
order deformation is given by \begin{equation*}
d^2=d^1-\phd{0,0,m-1}2m^2t^2\theta^3. \end{equation*} We compute
\begin{align*} \tfrac12[d^2,d^2]=&-\psd{1,0,0}2
t^1\theta^1-\psd{1,0,0}3 t^2\theta^2
+\phd{0,0,m}2(\theta^1\theta^3-m\theta^2\theta^3)\\
&-\psd{1,0,m}2(mt^3\theta^3+t^4\theta^1-mt^4\theta^2)
+\psd{1,0,m-1}2(m-1)m^2t^2t^3\theta^3 \\
&+\psd{1,0,m-2}2(m-1)m^2(t^2)^2\theta^3\\&
+\phd{0,0,m-1}2(-\theta^1+(m-1)\theta^2)m^2t^2\theta^3. \end{align*}
Notice that $\phd{0,0,m-1}2$ is not even a cocycle, so at first this
may seem to be a problem, because we know the right hand side is a
cocycle, and all the other terms are cocycles.  The answer is that by
using the second order relations we can see that the coefficient of
this term is zero.

Second, note that the original coboundary term $\psd{1,0,m-1}2$
reappears, as well as the new coboundary term $\psd{1,0,m-2}2$. In
fact, it is not hard to see that as we take higher order deformations
of $d_c$, eventually every coboundary of the form $\psd{1,0,k}2$ with
$k\le m-1$ will appear. It is not hard to guess that a miniversal
deformation $d^\infty$ is given by \begin{align*}
d^\infty=d^1+\sum_{k=0}^{m-1}\phd{1,0,k}2x^k, \end{align*} where $x^k$
is a power series in $t^3$, multiplied by $(t^2)^{m-k}\theta^3$. To see
this, suppose that $d^\infty$ has the form above, and compute
\begin{align*} \tfrac12[d^\infty,d^\infty]=&-\psd{1,0,0}2
t^1\theta^1-\psd{1,0,0}3 t^2\theta^2
+\phd{0,0,m}2(\theta^1\theta^3-m\theta^2\theta^3)
\\&-\psd{1,0,m}2(mt^3\theta^3+t^4\theta^1-mt^4\theta^2)+\sum_{k=1}^{m-1}\phd{1,0,k}2x^k(-\theta^1+k\theta^2)
\\&+\psd{1,0,m-1}2(-mt^2\theta^3+x^{m-1}(\tfrac1m-(m-1)t^3) \\
&+\sum_{k=1}^{m-2}\psd{1,0,k}2((\tfrac{m-k}m-kt^3)x^k-(k+1)x^{k+1}) .
\end{align*}
 From this computation one obtains that
 \begin{align*}
x^{m-1}=\frac{m^2t^2\theta^3}{1-m(m-1)t^3},\qquad x^k=\frac{(k+1)mt^2
x^{k+1}}{m-k-mkt^3},\quad k<m-1.\\ \end{align*} The second order
relations remain unmodified, and using them, we obtain that the term
involving $\phd{1,0,k}2$ is zero.  Thus we get a closed form
expression for the miniversal deformation,  even though it cannot be
obtained as a finite order deformation.

 \subsection{A miniversal
deformation of $\dt 0$} Since there are an infinite number of
cohomology classes, let us give them some simple labels.  Let
\begin{equation*} \xi=\psd{1,0,0}2,\quad\sigma=\phd{0,1,0}2\quad
\psi_k=\psd{1,0,k-1}3,\quad\phi_k=\phd{0,0,k}3\quad k\ge 1
\end{equation*} The universal infinitesimal deformation is given by
\begin{equation*} d^1=\xi s+\sigma \eta +\psi_kt^k+\phi_k\theta^k,
\end{equation*} where $s,t^k$ are even and $\eta,\theta^k$ are odd
parameters. Then \begin{align*} \tfrac12[d^1,d^1]&=-\xi
s\eta+\psi_{k+l-1}(k-l-1)t^k\theta^l+\tfrac12\phi_{k+l-1}(k-l)\theta^k\theta^l.
\end{align*} Thus $[d^1,d^1]$ has no coboundary terms, so the universal
infinitesimal deformation is miniversal, and the relations are given by
\begin{gather*} s\eta=0\\ \sum_{k+l=n+1}(k-l-1)t^k\theta^l=
\sum_{k+l=n+1}(k-l)\theta^k\theta^l=0,\qquad n=1\cdots\\ \end{gather*}

\subsection{A miniversal deformation of $d_c$ when $c$ is a negative
rational number} Decomposing $\frac c{c-1}=\frac rs$ as before, we
obtain an infinite number of cohomology classes, which we label as
follows. \begin{align*}
\xi_1=&\psd{1,0,0}2\quad&\xi_2=&\psd{1,0,0}3\quad&\sigma_1=&\phd{0,1,0}2\\
\psi_k=&\psd{1,kr,k(s-r)+1}3&\phi_k=&\phd{0,kr,k(s-r)+1}3\quad k\ge 0\\
\end{align*} The universal infinitesimal deformation is given by
\begin{equation*}
d^1=\xi_1s^1+\xi_2s^2+\sigma_1\eta^1+\psi_kt^k+\phi_k\theta^k,
\end{equation*} where $s^i, t^i$ are even and $\eta^1, \theta^k$ are
odd parameters. The bracket calculations needed to compute $[d^1,d^1]$
are \begin{align*}
[\xi_1,\sigma_1]=&-\xi_1\quad&[\xi_1,\phi_k]=&-\psd{1,kr-1,k(s-r)+1}3kr\\
[\xi_2,\sigma_1]=&0\quad&[\xi_2,\phi_k]=&-\psd{1,kr,k(s-r)}3(k(s-r)+1)\\
[\psi_k,\sigma_1]=&\psi_kkr\quad&[\psi_k,\phi_l]=&\psi_{k+l}(k-l)(s-r)\\
[\sigma_1,\phi_k]=&-\phi_kkr\quad&[\phi_k,\phi_l]=&\phi_{k+l}(k-l)(s-r)
\end{align*} The coboundary terms appearing above are \begin{align*}
-\psd{1,kr-1,k(s-r)+1}3kr=&D\left(\phd{0,kr-1,k(s-r)+1}3kr\right),\quad
k>0\\
-\psd{1,kr,k(s-r)}3(k(s-r)+1)=&D\left(\phd{0,kr,k(s-r)}3\tfrac{(k(s-r)+1)(s-r)}r\right),\quad
k>0. \end{align*} The second order relations are \begin{gather*}
s^1\eta^1=0,\qquad s^2\theta^0=0\\ nrt^n\eta^1+(s-r)\sum_{k+l=n}
t^k\theta^l(k-l)=0,\quad n=0\cdots\\
-2nr\eta^1\theta^n+(s-r)\sum_{k+l=n} \theta^k\theta^l(k-l)=0,\quad
n=0\cdots \end{gather*} The factor 2 appearing in the first summand of
the second relation occurs because the terms $[\phi_k,\phi_l]$ appear
only once in the bracket $[d^1,d^1]$ because they are like terms. (Note
that all the other like terms have zero self-brackets.) The second
order deformation is given by \begin{equation*}
d^2=d^1+\gamma_ks^1\theta^k+\epsilon_ks^2\theta^k, \end{equation*}
where $\gamma_k=\phd{0,kr-1,k(s-r)+1}3kr$ and
$\epsilon_k=\phd{0,kr,k(s-r)}3\tfrac{(k(s-r)+1)(s-r)}r$ for $k\ge 1$.
In order to compute $[d^2,d^2]$ we need to calculate the following
brackets. \begin{align*}
&[\xi_1,\gamma_k]=-\psd{1,kr-2,k(s-r)+1}3kr(kr-1)\\
&[\xi_2,\gamma_k]=-\psd{1,kr-1,k(s-r)}3kr(k(s-r)+1)\\
&[\sigma_1,\gamma_k]=-\phd{0,kr-1,k(s-r)+1}3kr(kr-1)\\
&[\psi_k,\gamma_l]=-\psd{1,(k+l)r-1,(k+l)(s-r)+1}3lr(l-k)(s-r)\\
&[\phi_k,\gamma_l]=-\phd{0,(k+l)r-1,(k+l)(s-r)+1}3lr(l-k)(s-r)\\
&[\gamma_k,\gamma_l]=\phd{0,(k+l)r-2,(k+l)(s-r)+1}3r^2lk(l-k)(s-r)\\
&[\xi_1,\epsilon_k]=-\psd{1,kr-1,k(s-r)}3k(s-r)(k(s-r)+1)\\
&[\xi_2,\epsilon_k]=-\psd{1,kr,k(s-r)-1}3\tfrac{k(s-r)^2(k(s-r)+1)}r\\
&[\sigma_1,\epsilon_k]=-\phd{0,kr,k(s-r)}3k(s-r)(k(s-r)+1)\\
&[\psi_k,\epsilon_l]=-\psd{1,(k+l)r,(k+l)(s-r)}3\tfrac{(l(s-r)+1)(s-r)((l-k)(s-r)-1)}r\\
&[\phi_k,\epsilon_l]=-\phd{0,(k+l)r,(k+l)(s-r)}3\tfrac{(l(s-r)+1)(s-r)((l-k)(s-r)-1)}r\\
&[\epsilon_k,\epsilon_l]=\phd{0,(k+l)r,(k+l)(s-r)-1}3\tfrac{(s-r)^3(k(s-r)+1)(l(s-r)+1)(l-k)}{r^2}\\
&[\gamma_k,\epsilon_l]=\phd{0,(k+l)r-1,(k+l)(s-r)}3k(s-r)(l(s-r)+1)((l-k)(s-r)-1)
\end{align*}
 Note that all the odd terms appearing above are
automatically cocycles, but none of the even terms are.  Therefore, the
sum of the even terms must be zero up to third order, using only the
third order relations. Since none of the terms appearing above are
cohomology classes, the third order relations are the same as the
second order ones.  Let us examine the even terms to see what is going
on.

The term in $[\sigma_1,\gamma_k]$ will be multiplied by the parameters $\eta^1$ and $s^1\theta^k$,
so it vanishes by the first relation $\eta^1s^1=0$. To see why the terms arising in $[\phi_k,\gamma_l]$
cancel, note that we need to sum the terms with $k+l=n$.  Adding the coefficients from the term
$[\phi_k,\gamma_l]$ with that from $[\gamma_k,\phi_l]$, we get the coefficient $nr(l-k)(s-r)s^1\theta^k\theta^l$.
Thus
\begin{equation*}
\sum_{k+l=n}lr(l-k)(s-r)s^1\theta^k\theta^l=\frac12\sum_{k+l=n}nr(l-k)(s-r)s^1\theta^k\theta^l.
\end{equation*}
If
you multiply the third relation by $s^1$, the first term drops out by the first relation, so it follows
that the sum above is zero.

The terms $[\gamma_k,\gamma_l]$, $[\epsilon_k,\epsilon_l]$ and $[\gamma_k,\epsilon_l]$
are multiplied by four parameters, so are
automatically zero up to third order.  Thus we are left with considering the terms $[\sigma_1,\epsilon_k]$
and $[\phi_k,\epsilon_l]$.  Note that these terms $[\sigma_1,\epsilon_n]$ and $[\phi_k,\epsilon_l]$ involve
the same cochain when $k+l=n$. Consider the coefficient arising from $[\phi_l,\epsilon_k]$, which is
$\frac{(k(s-r)+1)((l-k)(s-r)-1)s^2}r\theta^l\theta^k$. When we add this to $[\phi_k,\epsilon_l]$ we
obtain the coefficient
\begin{equation*}
-(k-l)(s-r)(n(s-r)+1)\theta^k\theta^l.
\end{equation*}
After multiplying this coefficient by $1/2$, summing, and adding the coefficient from $[\sigma,\epsilon_k]$, we obtain
exactly zero by the third relation.

While all these cancellations seem to appear miraculously, it really is not necessary to carry out this verification.
By the results in \cite{fp1}, the construction of the miniversal deformation which we are engaging in
is guaranteed to work,  so that any terms which arise in the brackets which are not cocycles must cancel.
Therefore, we can ignore these terms and concentrate on the coboundary terms above.

Of the five types of odd coboundary terms appearing in $[d^2,d^2]$,
$\psd{1,kr-2,k(s-r)+1}3$, $\psd{1,kr-1,k(s-r)}3$, and
$\psd{1,kr,k(s-r)-1}3$ are new types, while the other two are
coboundaries of $\gamma$ and $\epsilon$ terms.  The parameters for the
new types are $(s^1)^2\theta^k$, $s^1s^2\theta^k$ and
$(s^2)^2\theta^k$, and the powers of $s^1$ and $s^2$ which occur can be
read by looking at how much the middle and third upper indices have
decreased from the values $kr$ and $k(s-r)+1$ in the cohomology class
$\psd{1,kr,k(s-r)+1}3$.  Notice that this observation is also true for
$\gamma_k$ and $\epsilon_k$ as well. Putting this information
together,  we now construct a miniversal deformation.

Let
\begin{equation*}
\beta_{k,x,y}=\phd{0,kr-x,k(s-r)+1-y}3,\quad
 0\le x\le kr,\quad0\le y\le k(s-r)+1.
\end{equation*}
Note that if $x\ge r$, then for $x'=x-r$, $k'=k+1$ and $y'=y+s-r$,
we have $\beta_{k,x,y}=\beta_{k+1,x',y'}$, so we only need to consider the case $0\le x<r$.
Also, note that $\phi_k=\beta_{k,0,0}$.
We claim that a miniversal deformation $d^\infty$ of
$\dt c$ is given by
\begin{equation*}
d^\infty=\xi_1s^1+\xi_2s^2+\psi_kt^k+\sigma_1\eta^1+\phi_0\theta^0+
\beta_{k,x,y}u^{k,x,y}\theta^k,
\end{equation*}
where we restrict ourselves to the case $0\le x<r$, where $u^{k,x,y}$ is a power
series in the parameters, which we will determine by a recursive process. Since
$\beta_{k,0,0}=\phi_k$, we know that $u^{k,0,0}=\theta^k$.

Consider the brackets \begin{align*}
[\xi_1,\beta_{k,x,y}]=&-\psd{1,kr-x-1,k(s-r)+1-y}3(kr-x)\\
[\xi_2,\beta_{k,x,y}]=&-\psd{1,kr-x,k(s-r)-y}3(k(s-r)+1-y)=\\
[\psi_k,\beta_{l,x,y}]=&\psd{1,(k+l)r-x,(k+l)(s-r)+1-y}3((k-l)(s-r)+y)\\
[\sigma_1,\beta_{k,x,y}]=&-\beta_{k,x,y}(kr-x)\\
[\phi_0,\beta_{k,x,y}]=&-\beta_{k,x,y}(k(s-r)-y)\\
[\beta_{k,x,y},\beta_{l,u,v}]=&\beta_{k+l,x+u,y+v}((k-l)(s-r)-(y-v)).
\end{align*} Now \begin{align*}
-\psd{1,kr-x-1,k(s-r)+1-y}3=&D(\beta_{k,x+1,y})\tfrac{s-r}{(x+1)(s-r)-ry}\\
-\psd{1,kr-x,k(s-r)-y}3=&D(\beta_{k,x,y+1})\tfrac{s-r}{x(s-r)-r(y+1)}\\
-\psd{1,(k+l)r-x,(k+l)(s-r)+1-y}3=&D(\beta_{k+l,x,y})\tfrac{s-r}{x(s-r)-ry}\\
\end{align*} except, of course, when the denominators on the right hand
side vanish. Since $r$ and $s$ are relatively prime, $(x+1)(s-r)-ry=0$
only when $x+1=r$ and $y=s-r$, in which case
$\psd{1,kr-x-1,k(s-r)+1-y}3=\psi_{k-1}$. Since $x<r$, the denominator
of the second fraction never vanishes. The denominator of the third
fraction only vanishes when $x=y=0$, in which case, since
$\beta_{l,0,0}=\ph_l$, the bracket is just $\psi_{k+l}(k-l)(s-r)$ as
computed earlier. From this, we compute \begin{align*}
\tfrac12[d^\infty,d^\infty]&=
-\xi_1s^1\eta^1-\xi_2s^2\theta^0+\psi_kkrt^k\eta^1+\psi_kk(s-r)t^k\theta^0\\
&-D(\beta_{k,x,y}u^{k,x,y})+D(\beta_{k,x+1,y})\left(\tfrac{(kr-x)(s-r)}{(x+1)(s-r)-ry}\right)s^1u^{k,x,y}\\&
+D(\beta_{k,x,y+1})\left(\tfrac{(k(s-r)+1-y)(s-r)}{x(s-r)-r(y+1)}\right)s^2u^{k,x,y}\\
&-D(\beta_{k+l,x,y})\left(\tfrac{((k-l)(s-r)+y)(s-r)}{x(s-r)-ry}\right)t^ku^{l,x,y}\\&
-\beta_{k,x,y}(kr-x)\eta^1u^{k,x,y}
-\beta_{k,x,y}(k(s-r)-y)\theta^0u^{k,x,y}\\&
+\tfrac12\beta_{k+l,x+u,y+v}((k-l)(s-r)-(y-v)))u^{k,x,y}u^{l,u,v} ,
\end{align*} 
except for the cases when the denominators above vanish. When $x=r-1$
and $y=s-r$, then the first coboundary term above is replaced by
$-\psi_{k+1}((k-1)r+1)s^2u^{k,r-1,s-r}$, and when $x=y=0$ the second
coboundary term is replaced by $\psi_{k+l}(k-l)(s-r)t^k\theta^l$.

From the equation above,  we can easily determine the relations on the
base of the miniversal deformation. In fact, only one of the relations
is modified from the second order relations. The $\beta_{k,0,0}$ terms
are cohomology classes, so they give rise to relations, and the other
cohomology classes are immediately identifiable.  We obtain the
relations \begin{gather*} s^1\eta^1=0,\qquad s^2\theta^0=0\\
nrt^n\eta^1-(nr+1)s^1u^{n+1,r-1,s-r}+(s-r)\!\sum_{k+l=n}\!t^k\theta^l(k-l)=0,
n=0\cdots\\ -2nr\eta^1\theta^n+(s-r)\sum_{k+l=n}
\theta^k\theta^l(k-l)=0,\quad n=0\cdots \end{gather*}

Let us determine the coefficients $u^{k,x,y}$. We already know
$u^{n,0,0}$. The coefficients of each coboundary term must vanish.
First, when $x=0$, by summing up all the terms involving
$D(\beta_{n,0,y+1})$ we obtain \begin{align*}
0=&(n(s-r)+1-y)(s-r)s^2u^{n,0,y}\\&-\sum_{k=1}^{n-1}((2k-n)(s-r)+y+1)(s-r)t^ku^{n-k,0,y+1}\\&+
(r(y+1)-(n(s-r)+y+1)t^0)u^{n,0,y+1}, \end{align*} which determines
$u^{n,0,y+1}$ in terms of $u^{n,0,y}$ and $u^{k,0,y+1}$ for $k<n$. 
This solution is a power series in $t^0$, and is polynomial in the
parameters $s^2$, $t^k$ and $\theta^k$ for $0<k<n$. Secondly, summing
up all coefficients involving $D(\beta_{n,x+1,y})$ we obtain that
\begin{align*} 0=&(nr-x-1)s^1u^{n,x,y}+(n(s-r)-y)s^2u^{n,x+1,y-1}\\
&-\sum_{k=1}^{n-1}((k-l)(s-r)+y)t^ku^{l,x+1,y}\\
&-\left(\tfrac{(x+1)(s-r)-ry}{s-r}+(n(s-r)+y)t^0\right)u^{n,x+1,y},
\end{align*} which determines $u^{n,x+1,y}$ in terms of $u^{n,x,y}$,
$u^{n,x+1,y-1}$ and the values of $u^{l,x+1,y}$ for $l<n$.  Note that
we obtain the solution for $u^{n,x+1,y}$ as a power series in $t^0$,
and that it is polynomial in the parameters $s^1$, $s^2$, $t^k$, and
$\theta^k$, again for $1<k<n$.

We have provided the details of the recursive construction here in order to illustrate the method involved.
For subsequent examples we will not provide such complete details.

\subsection{A miniversal deformation of $\dstar$}
Since there are a lot of cohomology classes,  let us make some abbreviations for the cohomology
classes.  Let
\begin{equation*}
\psi_k=\psd{1,0,k-1}2\qquad
\xi_k=\psd{1,0,k-1}3\\
\end{equation*}
be the odd cohomology classes of degree $k$,
\begin{equation*}
\phi_k=\phd{0,0,k}3+\phd{1,0,k-1}1k\qquad
\sigma_k=\phd{0,k,0}2+\phd{0,k-1,1}3
\end{equation*}
be the even cohomology classes of degree $k$,
and $\tau=\phd{0,1,0}3$ be the extra even cohomology class in $L_1$.
Then the universal infinitesimal deformation of $\dstar$ is given by
\begin{equation*}
d^1=\dstar+\psi_ks^k+\xi_kt^k +\tau\zeta+\phi_k\theta^k+\sigma_k\eta^k,
\end{equation*}
where $s^k$ and $t^k$ are even parameters, and $\zeta$, $\theta^k$ and
$\eta^k$ are odd parameters.
Computing the brackets necessary to compute $[d^1,d^1]$, and comparing them to the coboundary calculations, we obtain
\begin{align*}
[\psi_k,\tau]=&\left(\psd{1,1,k-2}2(k-1)-\psd{1,0,k-1}3\right)
=\begin{cases}
D(\phd{0,0,k-1}2)& k>1\\
-\xi_1&k=1
\end{cases}\\
[\psi_k,\phi_l]=&\psi_{k+l-1}(k+l-1)\\
[\psi_k,\sigma_l]=&-\psd{1,l-1,k-1}2(l+1-k)-\psd{1,l-2,k}3(l-1)\\
=&\begin{cases}
\psi_k(k-2)& l=1\\
D\left(-\phd{0,0,k}2\left(\tfrac{3-k}k\right)\right)+\psi_{k+1}\tfrac{3-2k}k &l=2\\
D\left(-\phd{0,l-3,k+1}3\left(\tfrac{kl-2k+l+1}{k(k+1)}\right)-\phd{0,l-2,k}2\left(\tfrac{l+1-k}k\right)\right)
&\text{otherwise}
\end{cases}\\
[\xi_k,\tau]=&\psd{1,1,k-2}3(k-1)=
\begin{cases}
D(\phd{0,0,k-1}3)& k>1\\
0&k=1
\end{cases}\\
[\xi_k,\phi_l]=&\xi_{k+l-1}(k-1)\\
[\xi_k,\sigma_l]=&\psd{1,l-1,k-1}3(k-2)=
\begin{cases}
D(\phd{0,l-2,k}3)\tfrac{k-2}k&l>1\\
\xi_k(k-2)&l=1
\end{cases}\\
[\tau,\phi_k]=&-(\phd{1,1,k-2}1(k-1)+\phd{0,1,k-1}3))k=
\begin{cases}
D(-\psd{0,0,k-1}1)k&k>1\\
-\tau&k=1
\end{cases}\\
[\tau,\sigma_k]=&0\\
[\phi_k,\phi_l]=&\phi_{k+l-1}(k-l)\\
[\phi_k,\sigma_l]=&\left(\phd{1,l-1,k-1}1k+\phd{0,l-1,k}3\right)(k-1)=\begin{cases}
D(\psd{0,l-2,k}1)(k-1)&l>1\\
\phi_k(k-1)&l=1
\end{cases}\\
[\sigma_k,\sigma_l]=&\sigma_{k+l-1}(k-l)
\end{align*}
The second order relations are
\begin{align*}
&s^1\zeta+t^1\eta^1=0\\
&\zeta\theta^1=0\\
&-s^1\eta^1+s^1\theta^1=0\\
&(n-2)s^n\eta^1+\left(\tfrac{5-2n}{n-1}\right)t^{n-1}\eta^2+n\sum_{k+l=n+1}s^k\theta^l=0,\qquad n> 1\\
&(n-2)t^n\eta^1+\sum_{k+l=n+1}(k-1)t^k\theta^l=0,\qquad n>1\\
&(n-1)\theta^n\eta^1+\sum_{k+l=n+1}(k-l)\theta^k\theta^l=0,n\ge1\\
&\sum_{k+l=n}(k-l)\eta^k\eta^l=0.
\end{align*}
The second order deformation is given by
\begin{align*}
d^2=&d^1
-\phd{0,0,k-1}2s^k\zeta
+\left(-\phd{0,l-3,k+1}3\left(\tfrac{kl-2k+l+1}{k(k+1)}\right)-\phd{0,l-2,k}2\left(\tfrac{l+1-k}k\right)\right)s^k\eta^l\\
&-\phd{0,0,k-1}3t^k\zeta
+\phd{0,l-2,k}3\tfrac{2-k}kt^k\eta^l
+\phd{0,0,k-1}1\zeta\theta^k
+\psd{0,l-2,k}1(1-k)\theta^k\eta^l
\end{align*}
In order to determine $[d^2,d^2]$, it is necessary to compute 51 additional brackets.  Many of the
brackets give coboundary terms, so it is clear that $d^2$ is far from the end of the
story.  The calculations do not seem especially illuminating, so we felt that it would not
be useful to provide them here.  It is not that surprising that the miniversal deformation
of $\dstar$ is a complicated object, given the complexity of its cohomology.

\section{Codifferentials of degree 2 of the second kind}

A codifferential of degree 2 of the second kind is of the form
$$d=\psd{0,2,0}1a+ \psd{0,1,1}1b + \psd{0,0,2}1c,$$ which we will say
is of type $(a,b,c)$. If either $a$ or $b$ is nonzero, then it is
clearly equivalent to a codifferential of type $(1,b',c')$, for some
$b'$ and $c'$.  Note that the only type which cannot be reduced in this
way is type $(0,b,0)$, which is clearly also of type $(0,1,0)$.
However, let us examine type $(0,1,0)$ to see what it is equivalent to.
Applying a standard linear automorphism, we obtain that
$d=\psd{0,1,1}1$ is equivalent to any codifferential of the form
$$d'=\psd{0,2,0}1 \frac{2rs}q +\psd{0,1,1}1\frac{ru+ts}q
+\psd{0,0,2}1\frac{2tu}q,\qquad ru-ts\ne 0, q\ne0.$$ If we set $q=2rs$,
$x=u/s$ $y=s/r$, then we obtain type $(1,b,c)$ where $2b=x+y$, $c=xy$,
and we must avoid the condition $xy=1$.  But this occurs exactly when
$b^2=c$. Thus type  $(0,1,0)$ is equivalent to type $(1,b,c)$ whenever
$b^2\ne c$.  It is also clear that type $(0,1,0)$ is not equivalent to
type $(1,0,0)$.

Let us next study type $(1,0,0)$.  Applying a linear automorphism, we
see that $\psd{0,2,0}1$ is equivalent to codifferentials of the form
$$
d'=\psd{0,2,0}1\frac {r^2}q+\psd{0,1,1}1\frac{rt}q+\psd{0,0,2}1\frac{t^2}q.
$$
If we set $r^2=q$, then we see that this $d$ is equivalent to any codifferential of the form $(1,b,b^2)$,
exactly the types not covered by the first case.

Thus there are only two types of codifferentials, represented by $\psd{0,2,0}1$ and $\psd{0,1,1}1$.
Let us study the second type first.
\subsection{Type $(0,1,0)$}
Let $D(\ph)=\br{\ph}{\psd{0,1,1}1}$. Then we obtain the following table of coboundaries.
\begin{align*}
D(\phd{1, q, n-q-1}1) &=\psd{0, 1+q, n-q}1\\
D(\phd{0, p, n-p}2) &=-\psd{0,p, n-p+1}1\\
D(\phd{0, p, n-p}3) &=-\psd{0,p+1, n-p}1\\
D(\psd{0, p, n-p}1) &=0\\
D(\psd{1, q, n-q-1}2) &=\phd{1,q,n-q}1+\phd{0,q+1,n-q}2\\
D(\psd{1, q, n-q-1}3) &=\phd{1,q+1,n-q-1}1+\phd{0,q+1,n-q}3\\
\end{align*}
Then we have $n+1$ odd cocycles of the form $\psd{0,p,n-p}1$,  and $2n$
even cocycles of the form $\phd{1,q,n-q}1+\phd{0,q,n-q-1}3$ and
$\phd{1,q,n-q-1}1+\phd{0,q+1,n-q-1}2$ This means $z_n=2n|n+1$, so
$b_{n}=2n|n+2$, and $h_n=2|0$ if $n>1$. We have
\begin{align*}
&H^1=\langle \psd{0,0,1}1,\psd{0,1,0}1,\phd{1,0,0}1+\phd{0,0,1}3,\phd{1,0,0}1+\phd{0,1,0}2\rangle\\
&H^n=\langle\phd{1,0,n-1}1+\phd{0,0,n}3,\phd{1,n-1,0}1+\phd{0,n,0}2\rangle,\qquad\text{if $n>1$}
\end{align*}
and all cohomology for $n>1$ is even. Thus we don't obtain any
deformations in the Lie algebra direction.
\subsection{Type $(1,0,0)$}
Let $D(\ph)=\br{\ph}{\psd{0,2,0}1}$. Then we obtain the following table
of coboundaries.
\begin{align*}
D(\phd{1, q, n-q-1}1) &=\psd{0, 2+q, n-q-1}1\\
D(\phd{0, p, n-p}2) &=-2\psd{0,p+1, n-p}1\\
D(\phd{0, p, n-p}3) &=0\\
D(\psd{0, p, n-p}1) &=0\\
D(\psd{1, q, n-q-1}2) &=2\phd{1,q+1,n-q-1}1+\phd{0,q+2,n-q-1}2\\
D(\psd{1, q, n-q-1}3) &=\phd{0,q+2,n-q-1}3\\
\end{align*}

Besides the obvious $n+1$ odd cocycles $\psd{0,p.n-p}1$ and  $n+1$ even
ones $\phd{0,p,n-p}3$, we also have $n$ more even cocycles
$$2\phd{1,q,n-q-1}1+\ph_2^{0,q+1,n-q-1},$$ so $z_n=2n+1|n+1$ and
$h_n=3|1$, if $n>1$.  In fact, it is easily seen that
\begin{align*}
&H^1=\langle \psd{0,0,1}1,\psd{0,1,0}1,\phd{0,0,1}3,\phd{0,1,0}3,2\phd{1,0,0}1+\phd{0,1,0}2\rangle\\
&H^n=\langle \psd{0,0,n}1,\phd{0,0,n}3,\phd{0,1,n-1}3,2\phd{1,0,n-1}1+\phd{0,1,n-1}2,\rangle,\quad\text{if $n>1$}
\end{align*}

Let us think about the moduli space of two points given by these
codifferentials.  Note that type $(1,0,\epsilon)$ is the same as type
$(0,1,0)$ for any nonzero value of $\epsilon$.  Thus type $(1,0,0)$ is
infinitesimally close to $(0,1,0)$, but not the other way around. It is
not surprising, therefore to see that type $(1,0,0)$ has a nontrivial
deformation as a Lie algebra.

\section{Miniversal deformations of degree 2 codifferentials of the
second kind}

\subsection{Miniversal deformations of Type (0,1,0)}
Let $d=\psd{0,1,1}1$, and let us label the cohomology classes as
follows:
\begin{equation*}
\psi_1=\psd{0,0,1}1,\ \ \psi_2=\psd{0,1,0}1,\
\phi_n=\psd{1,0,n-1}1+\psd{0,0,n}3,\quad \sigma_n=\phd{1,n-1,0}1+\phd{0,n,0}2
\end{equation*}
The universal infinitesimal deformation is given by
\begin{equation*}
d^1=\psd{0,1,1}1+\psi_1 s^1 +\psi_2s^2+\phi_n \theta^n +\sigma_n\eta^n,
\end{equation*}
where $s^i$ are even parameters and $\theta^n$ and $\eta^n$ are odd parameters.
The brackets we need to compute are
\begin{align*}
[\psi_1,\phi_k]=&[\psi_2,\sigma_k]=[\phi_k,\sigma_l]=0\\
[\psi_1,\sigma_k]=&-\psd{0,k-1,1}1=
\begin{cases}
D(\phd{0,k-1,0}2)&k>1\\
-\psi_1&k=1
\end{cases}\\
[\psi_2,\phi_k]=&-\psd{0,1,k-1}1=
\begin{cases}
D(\phd{0,1,k-2}2)& k>1\\
-\psi_2&k=1
\end{cases}\\
[\phi_k,\phi_l]=&\phi_{k+l-1}(k-l),\quad
[\sigma_k,\sigma_l]=\sigma_{k+l-1}(k-l)\\
\end{align*}
Thus the second order relations are
\begin{gather*}
s^1\eta^1=s^2\theta^1=0\\
\sum_{k+l=n+1}(k-l)\theta^k\theta^l=0,\quad n\ge 1\\
\sum_{k+l=n+1}(k-l)\eta^k\eta^l=0,\quad n\ge 1\\
\end{gather*}
Let $\gamma_k=-\phd{0,k-1,0}2$ and $\epsilon_k=-\phd{0,1,k-2}2$.
Then the second order deformation is
\begin{equation*}
d^2=d^1+\gamma_k s^1\eta^k+\epsilon_k s^2\theta^k, k>1.
\end{equation*}
The brackets needed to compute $[d^2,d^2]$ are
\begin{align*}
[\psi_1,\gamma_k]=&[\phi_k,\gamma_k]=[\psi_1,\epsilon_k]=[\epsilon_k,\epsilon_l]=0\\
[\psi_2,\gamma_k]=&-\psd{0,k-1,0}1\qquad\qquad\qquad[\psi_2,\epsilon_k]=-\psd{0,1,k-1}1\\
[\gamma_k,\gamma_l]=&\phd{0,k+l-3,0}2(k-l)\qquad\quad
[\phi_k,\epsilon_l]=\phd{0,1,k+l-2}2(l-1)\\
[\sigma_k,\gamma_l]=&\phd{1,k+l-3,0}1(1-k)+\phd{0,k+l-2,0}2(l-k-1)\\
[\sigma_k,\epsilon_l]=&\phd{1,k-1,l-1}1(1-k)+\phd{0,k,l-1}2(1-k)
\end{align*}
With the exception of the terms $[\psi_2,\gamma_2]$ and $[\sigma_k,\gamma_l]$, the terms above do not
involve any cohomology classes. Most of the terms vanish, after taking into account the
third order relations, which are the same as the second order relations except that the relation
$s^2\theta^1=0$ becomes $s^2\theta^1+s^1s^2\eta^1=0$.  The modification of the second order relations
by addition of higher order terms is a common pattern that occurs in the construction of miniversal
deformations, as was illustrated by several examples in \cite{fp2}.

The coboundary terms which do not vanish are \begin{align*}
-&\psd{0,k-1,0}1=D\left(\phd{0,k-2,0}3\right), k>2\\
-&\psd{0,1,k-1}1=D\left(\phd{0,0,k-1}3\right)\\
&\phd{1,k-1,l-1}1(1-k)+\phd{0,k,l-1}2(1-k)=D\left(\psd{1,k-1,l-2}2(1-k)\right).
\end{align*} Thus we have \begin{equation*}
d^3=d^2-\phd{0,k-2,0}3s^1s^2\eta^k-\phd{0,0,k-1}3 (s^2)^2\theta^k
-\psd{1,k-1,l-2}2(1-k)s^2\eta^k\theta^l. \end{equation*}
 The next step
would involve calculating brackets of these three new terms with all
the terms introduced in $d^1$ and $d^2$, as well as with each other. We
did not work out the construction here. However, in the next example,  we will
give some more detail,  showing how the process can be reduced to a
recursion.
 \subsection{Miniversal deformations of Type (1,0,0)}
 Let
$d=\psd{0,2,0}1$, and let us label the cohomology classes as follows
\begin{gather*} \xi=\psd{0,1,0}1\\ \psi_n=\psd{0,0,n}1,\quad
\phi_n=\phd{0,0,n}3,n>0\\
\sigma_n=\phd{0,1,n-1}3,\quad\tau_n=2\phd{1,0,n-1}1+\phd{0,1,n-1}2,\quad
n> 0 \end{gather*} The universal infinitesimal deformation is given by
\begin{equation*} d^1=\psd{0,2,0}1+\xi s^1 +\psi_n t^n+\phi_n \theta^n
+\sigma_n\eta^n+\tau_n\zeta^n, \end{equation*} where $s^1$ and $t^n$
are even parameters and $\theta^n$, $\eta^n$ and $\zeta^n$ are odd
parameters. The brackets we need to compute are \begin{align*}
[\xi,\phi_k]=&[\xi,\sigma_k]=[\tau_k,\tau_l]=0\\
[\xi,\tau_k]=&-\psd{0,1,k-1}1,\qquad [\psi_k,\phi_l]=\psi_{k+l-1}k\\
[\psi_k,\sigma_l]=&\psd{0,1,k+l-2}1k,\qquad
[\psi_k,\tau_l]=-2\psi_{k+l-1}\\
[\phi_k,\phi_l]=&\phi_{k+l-1}(k-l)\quad
[\phi_k,\sigma_l]=\sigma_{k+l-1}(k-l+1)\\
[\phi_k,\tau_l]=&\tau_{k+l-1}(1-l),\quad
[\sigma_k,\sigma_l]=\phd{0,2,k+l-3}3(k-l)\\
[\sigma_k,\tau_l]=&2\phd{1,1,k+l-3}1(1-l)+\phd{0,2,k+l-3}2(1-l)+\phd{0,1,k+l-2}3\\
\end{align*} Note that there are a few special cases: \begin{align*}
[\xi,\tau_1]=-\xi,\quad [\psi_1,\sigma_1]=\xi,\quad
[\sigma_1,\tau_1]=\sigma_1\\ \end{align*} Thus the second order
relations are \begin{gather*} -s^1\zeta^1+t^1\eta^1=
\theta^1\eta^1+\eta^1\zeta^1=0\\
\sum_{k+l=n+1}(k-l-1)\theta^k\eta^l+\eta\zeta=0,\quad n=2\cdots\\
\sum_{k+l=n+1}(k-l)\theta^k\theta^l=0,n=1\cdots\\
\sum_{k+l=n+1}(1-l)\theta^k\zeta^l= \sum_{k+l=n+1}k t^k
(\theta^l-2\eta^l)=0,n=1\cdots\\ \end{gather*} In order to set up the
recursion relations,  let \begin{equation*}
\gamma_{k,l}=\phd{0,k,l}2,\quad \alpha_{k,l}=\psd{1,k,l}2,\quad
\beta_{k,l}=\psd{1,k,l}3 \end{equation*} The coboundary terms arising
in $[d^1,d^1]$ can all be expressed in terms of coboundaries of these
cochains. \begin{align*} [\xi,\tau_k]=-\tfrac12D(\gamma_{0,k-1}),\quad
k>1\\ [\psi_k,\sigma_l]=\tfrac12D(\gamma_{0,k+l-2}),\quad k+l>2\\
[\sigma_k,\sigma_l]=D(\beta_{0,k+l-3})(k-l), \quad k+l>2\\
[\sigma_k,\tau_l]=D(\alpha_{0,k+l-3}(1-l))+\sigma_{k+l-1}\\
\end{align*}
 It is also easy to express the brackets of the $\gamma$,
$\beta$ and $\alpha$ cochains with each other and all of the cohomology
classes.  The resulting cochains can either be expressed in terms of
these cochains and cocycles, which makes it possible to find the
recursion relations.  We do not include the details of these
calculations here.
 \section{Conclusions}
In this paper,  we have explored the construction of all linear and
quadratic $\Z_2$-graded \linf\ structures on a $2|1$-dimensional space.  We also
provided some information about \linf\ structures in general on this
space.  Actually, we showed that there are two kinds of
codifferentials, those of the first and second kinds , and that any
codifferential is either a sum of cochains of the first kind, or a sum
of cochains of the second kind, and every such sum is a codifferential.
In this sense,  we have described all \linf\ structures on a $2|1$
dimensional space.  However, because we did not address equivalence in
general, we are a long way from classification of the \linf\
structures.

For linear and quadratic \linf\ structures, we did give a complete
classification, as well as classifying all \linf\ structures with a
leading linear term (they are equivalent to the structure given by the
leading term). For \linf\ structures with a leading quadratic term,  it
is necessary to classify the extensions of the quadratic
codifferentials, and that task we leave to a separate paper.  We have
not studied the problem of classifying \linf\ structures with a leading
term of degree 3 or higher, or even classified the structures of a
fixed degree larger than 2.

In addition,  we have shown how to construct a miniversal deformation for each of the equivalence classes
of degree 2 codifferentials,  providing complete details when the cohomology is not too complicated, and
indicating some methods of computation in general.  The main goal of this paper is not to provide an exhaustive
description of a miniversal deformation, but rather to give the reader an idea of the process involved
in its computation.

\providecommand{\bysame}{\leavevmode\hbox to3em{\hrulefill}\thinspace}
\providecommand{\MR}{\relax\ifhmode\unskip\space\fi MR }
\providecommand{\MRhref}[2]{%
  \href{http://www.ams.org/mathscinet-getitem?mr=#1}{#2}
}
\providecommand{\href}[2]{#2}

\end{document}